\pdfoutput=1
\documentclass[a4paper,11pt]{article}
\usepackage[english]{babel}
\usepackage{amsfonts}
\usepackage{amsthm}
\usepackage{amsmath}
\usepackage{amssymb}
\usepackage{enumerate}
\usepackage{mathtools}
\usepackage{tikz}
\usepackage{tikz-cd}
\usepackage[all]{xy}
\usepackage{graphicx}
\usepackage{amsmath}
\usepackage{adjustbox}
\usepackage{graphicx}
\usepackage{float}
\usepackage[toc]{appendix}
\usepackage[left=3.1cm,top=2.5cm,right=3.1cm,bottom=2.5cm]{geometry}
\usepackage{hyperref}
\usepackage{resizegather}
\usepackage[activate={true, nocompatibility}, final, tracking=true, kerning=true, spacing=true, factor=1100, stretch=10, shrink=10]{microtype}
\usepackage{dsfont}
\usetikzlibrary{matrix}
\usepackage{cleveref}
\usepackage{ulem}

\usepackage{quiver} 

\usepackage{hyperref}
\hypersetup{
colorlinks=true,
linkcolor=[RGB]{0,0,204},
filecolor=magenta,      
urlcolor=cyan,
citecolor=[RGB]{0,204,0},
}

\theoremstyle{plain}
\newtheorem{thm}{Theorem}[section]
\newtheorem{coroll}[thm]{Corollary}
\newtheorem{defn}[thm]{Definition}
\newtheorem{lemma}[thm]{Lemma}
\newtheorem*{claim}{Claim}
\newtheorem{fact}[thm]{Fact}
\newtheorem{example}[thm]{Example}
\newtheorem{prop}[thm]{Proposition}

\newtheorem{remark}[thm]{Remark}

\newtheorem*{mainthm*}{Main Theorem}
\newtheorem{notn}[thm]{Notation}

\newtheorem*{coroll*}{Corollary}
\newtheorem*{thm*}{Theorem}
\newtheorem{property}[thm]{Property}

\tikzset{
  symbol/.style={
    draw=none,
    every to/.append style={
      edge node={node [sloped, allow upside down, auto=false]{$#1$}}}
  }
}

\microtypecontext{spacing=nonfrench}

\DeclareMathOperator{\Sym}{Sym}

\DeclareMathOperator{\QCoh}{QCoh}

\DeclareMathOperator{\Gr}{Gr}
\DeclareMathOperator{\MHM}{MHM}

\newcommand{\bC}{\mathbb{C}}

\newcommand{\bQ}{\mathbb{Q}}

\newcommand{\bZ}{\mathbb{Z}}

\newcommand{\on}{\operatorname}

\newcommand{\Coh}{\on{Coh}}

\newcommand{\Ext}{ \on{Ext}}

\newcommand{\Hom}{ \on{Hom}}

\newcommand{\Spec}{\on{Spec}}

\newcommand{\Bun}{ \on{Bun} } 
\newcommand{\Pic}{\on{Pic}}

        \newcommand{\cC}{{\mathcal C}}
        
        \newcommand{\cE}{{\mathcal E}}
        \newcommand{\cF}{{\mathcal F}}
        
        \newcommand{\cH}{{\mathcal H}}
        \newcommand{\cI}{{\mathcal I}}

        \newcommand{\cL}{{\mathcal L}}
        \newcommand{\cM}{{\mathcal M}}
        \newcommand{\cN}{{\mathcal N}}
        \newcommand{\cO}{{\mathcal O}}

        \newcommand{\cS}{{\mathcal S}}
        \newcommand{\cT}{{\mathcal T}}
        \newcommand{\cU}{{\mathcal U}}
        \newcommand{\cV}{{\mathcal V}}
        
        \newcommand{\cX}{{\mathcal X}}
        \newcommand{\cY}{{\mathcal Y}}
        \newcommand{\cZ}{{\mathcal Z}}

        \newcommand{\bL}{{\mathbb L}}
        \newcommand{\git}{/\!\!/}

\usepackage{color} 
\definecolor{darkgreen}{rgb}{0.0, 0.7, 0.0}

\begin{document}
\title{\textbf{Relative \'etale slices and cohomology of moduli spaces}}
\author{Mark Andrea de Cataldo, Andres Fernandez Herrero and Andrés Ibá\~nez Nú\~nez}
\date{}
\maketitle
\begin{abstract}
We use techniques of Alper--Hall--Rydh to prove a local structure theorem for smooth morphisms between smooth stacks around points with linearly reductive stabilizers. This implies that the good moduli space of a smooth stack over a base has equisingular fibers. As an application, we show that any two fibers have isomorphic $\ell$-adic cohomology rings and intersection cohomology groups. If we work over the complex numbers, we show that the family is topologically locally trivial on the base, and that the intersection cohomology groups of the fibers fit into a polarizable variation of pure Hodge structures. We apply these results to derive some consequences for the moduli spaces of $G$-bundles on smooth projective curves, and for certain moduli spaces of sheaves on del Pezzo surfaces.
\end{abstract}

\tableofcontents

 \section{Introduction}

 The cohomology of smooth projective varieties is locally constant in families. More precisely, if $X \to S$ is a smooth projective family of varieties over a connected base $S$, then the ($\ell$-adic) cohomology groups of any two fibers are non-canonically isomorphic via specialization. If we work over the complex numbers, then the singular cohomology groups of the family $X \to S$ fit into a polarizable variation of pure Hodge structures. In particular, certain invariants, such as the Hodge numbers, are locally constant in smooth projective families.

 In moduli theory, it is common to encounter moduli problems that naturally vary over a base $S$. Such moduli problems are often represented by an algebraic stack $\cX \to S$. One can then employ the techniques of Geometric Invariant Theory \cite{seshadri-relative-git}, or more modern stack-theoretic methods \cite{_Alper_Existenceofmodulispacesforalgebraicstacks}, to construct a good moduli space $M \to S$ for the stack $\cX$. Even when the stack $\cX$ is smooth over $S$, the moduli space $M \to S$ might have singular fibers. In practice, this arises when there are strictly semistable objects in our moduli problem.
 
 The main observation in this paper is that the desirable behavior of the cohomology of smooth projective families also holds for projective morphisms $M \to S$ that are obtained as good moduli spaces of smooth stacks $\cX \to S$.

 To that effect, we first prove a local structure theorem for smooth morphisms of smooth stacks at points with linearly reductive stabilizers (\Cref{proposition: relative smooth etale slice}). It states that the \'etale local tubular neighborhoods of smooth points of stacks obtained in \cite[Thm. 1.2]{alper-hall-rydh-etale-slice} are functorial with respect to smooth morphisms. \Cref{proposition: relative smooth etale slice} should be thought of as a relative version of Luna \'etale slice theorem in the spirit of \cite{alper-hall-rydh-theetalelocalstructureofalgebraicstacks}. Alper--Janda have given an alternative proof of \Cref{proposition: relative smooth etale slice} \cite[Corollary 5.1]{jarod-felix-paper}.

 As a consequence of \Cref{proposition: relative smooth etale slice}, we obtain the following (see \Cref{cor: equisingularity} for a more precise formulation).
 \begin{thm} \label{thm: A intro} Let $k$ be an algebraically closed field. Let $\cX \to S$ be a smooth stack with affine relative diagonal over a smooth algebraic space $S$ over $k$. If $\cX$ admits a good moduli space $M$, then the fibers of $M \to S$ are equisingular.
 \end{thm}

 This equisingularity result, combined with the machinery of vanishing cycles, allows us to relate the cohomology of such moduli spaces as they vary in a family. For instance, we obtain the following simple description of the decomposition theorem, generalizing the classical theorem of Deligne on decomposition for projective smooth morphisms \cite{deligne-suites-spectrales}. 

 \begin{thm}[\Cref{thm: coh of gms in families}] \label{thm: B intro}
     Let $f: M \to S$ be a proper morphism of schemes such that $S$ is smooth over an algebraically closed field. Suppose that $M$ is the good moduli space of a smooth algebraic stack $\cX \to S$ with affine relative diagonal.
     Then the pushforward $Rf_* \mathcal{IC}_M$ of the topologist's intersection complex $\mathcal{IC}_{M}$ (as in Subsection \ref{preliminaries l-adic}) admits a direct sum decomposition 
     \[
Rf_* \mathcal{IC}_M \simeq \bigoplus_{i\geq 0}
(R^if_* \mathcal{IC}_M) [-i],
\]
and the higher pushforward sheaves  $R^if_* \mathcal{IC}_M$ are lisse sheaves on $S$ for all $i$.
 \end{thm}

\Cref{thm: B intro} is stated in the context of $\ell$-adic cohomology with coefficients in $\overline{\mathbb{Q}}_{\ell}$. It implies in particular that the $\ell$-adic cohomology (resp. intersection cohomology) groups of any two geometric fibers of $M \to S$ are non-canonically isomorphic if $S$ is connected; see \Cref{iso coh fibers}.

 If the ground field is $\mathbb{C}$, then the same results hold for singular cohomology (resp. intersection cohomology) with rational coefficients for the classical topology.

 \begin{thm}[\Cref{cor: equisingularity}; Propositions \ref{prop: mixed hodge module prop} and \ref{loc triv}; \Cref{thm: two fibres of locally isotrivial family are cohomologically equivalent singular cohomology}] \label{introthm: C}
     Let $f: M \to S$ be a projective morphism of $\mathbb{C}$-varieties such that $S$ is smooth and connected. Suppose that $M$ is the good moduli space of a smooth algebraic stack $\cX \to S$ with affine relative diagonal. Then, the following hold; 
     \begin{enumerate}[(a)]
         \item The family $f:M \to S$ is topologically locally trivial. In particular, all $S$-fibers are homeomorphic.
         \item The pushforward $Rf_* \mathcal{IC}_M$ of the topologist's intersection complex $\mathcal{IC}_{M}$ with rational coefficients admits a direct sum decomposition in the category of mixed Hodge modules
     \[
Rf_* \mathcal{IC}_M \simeq \bigoplus_{i\geq 0}
(R^if_* \mathcal{IC}_M) [-i].
\]
Every higher pushforward sheaf $R^if_* \mathcal{IC}_M$ is a local system underlying a polarizable variation of Hodge structures. The Hodge numbers $h^{p,q}(I\!H^{\bullet})$ associated with the intersection cohomology of any two $S$-fibers agree.
\item The higher pushforwards $R^if_*\mathbb{Q}_M$ of the constant sheaf $\mathbb{Q}_M$ are local systems underlying an admissible variation of mixed Hodge structures. The mixed Hodge numbers $h^{p,q}\Gr^k(H^{\bullet})$ associated with the singular cohomology of any two $S$-fibers agree.
     \end{enumerate}
 \end{thm}

In \Cref{sec: applications},
we derive three 
cohomological consequences as applications of our results.

 For our first application, in Subsection \ref{section: application G-bundles} we fix a complex reductive group $G$ and an element of the algebraic fundamental group $d \in \pi_1(G)$. We show that, for any two complex smooth projective curves $C_1, C_2$ of the same genus, the corresponding moduli spaces $M^d_G(C_1), M^d_G(C_2)$ of semistable $G$-bundles of degree $d$ have the same mixed Hodge numbers (resp. Hodge numbers for intersection cohomology). The same holds for the moduli of parabolic $G$-bundles on pointed curves of the same genus for any fixed choice of parabolic weights.

 In Subsection \ref{application: negatively polarized surfaces} we show a similar result for moduli of Gieseker semistable torsion-free sheaves on smooth surfaces $X$, where semistability is taken with respect to a polarization $\cL$ satisfying $\cL \cdot \omega_X<0$ (we call these $\omega$-negatively polarized surfaces). As we vary the surface and polarization on connected families, the cohomological and Hodge-theoretic invariants of the corresponding moduli spaces of sheaves do not vary, even in the case when the stability parameters are not generic (i.e. in the presence of strictly semistable sheaves). In the case of sheaves of pure dimension 1, we show that the induced isomorphisms on intersection cohomology are compatible with the perverse filtrations induced by the Fitting support morphism (also known as the Le Potier morphism).

 \subsection*{Notation}
 
We work over a fixed algebraically closed ground field $k$. All of the stacks in this paper are algebraic stacks over $k$, and all morphisms are $k$-morphisms. 
Given stacks $\cX, \cY$ over $k$, we write $\cX \times \cY$ to denote the fiber product over $\Spec(k)$. For an algebraic space $X$ endowed with an action of an algebraic group $G$, we write $X/G$ for the quotient stack (which is denoted by $[X/G]$ in works of other authors). We use cohomological grading conventions for all triangulated categories. If $\cX$ is a stack of finite type over $k$, by its derived category of quasicoherent sheaves we mean the derived category of the Grothendieck Abelian category of quasicoherent sheaves on $\cX.$
In fact, every time that we employ the quasicoherent derived category,
\cite[Cor. 9.2(5)]{hall-rydh-perfect-complexes-stacks} applies, so that
 there should  be no ambiguity with other versions of the derived category.

\subsection*{Acknowledgements} 
We would like to thank the referee for very helpful comments on how to improve the structure of the paper.
 We would like to thank Jarod Alper, Nathan Chen, Daniel Halpern-Leistner, Jochen Heinloth, Tudor Padurariu, David Rydh and Yukinobu Toda for helpful comments and discussions. The first-named author  has been supported  by a grant from the Simons Foundation (672936, de Cataldo)
 and by NSF Grants DMS 1901975 and DMS 2200492. The third-named author would like to express gratitude to the Mathematical Institute at the University of Oxford for its support.
 
 \section{Relative \'etale slices}

\subsection{Preliminaries and notation}

Let $\cX$ be an algebraic stack, quasi-separated and locally of finite type  over $k$. Let $x \in \mathcal{X}(k)$ be a closed point in $\cX$, and let $i\colon BG_x\to \cX$ be the residual gerbe at $x$, where $G_x$ is the stabilizer group of $x$. The closed immersion $i$ is cut out by an ideal $\cI \subset \cO_{\cX}$. The \textit{normal space}  $N_x$
to $x$ is the $k$-vector space $ (\cI/\cI^2)^\vee = 
(i^*\cI)^{\vee}$ with the structure of $G_x$-representation. If 
$x$ is a smooth point and $G_x$ is smooth over $k,$ then $N_x$ 
is the tangent space $T_{\cX, x}$ of the stack at $x$. We 
denote $\cN_x:=N_x/G_x=\Spec_{BG_x}(\Sym(i^* \cI))$ and call it 
the \textit{normal stack} of $\cX$ at $x$. Let $0 \in \cN_x(k)$
denote the point corresponding to the origin $0$ in the vector 
space $N_x$.
If $\cY$ is another algebraic stack which is quasi-separated and locally of finite type over $k$, and $f\colon \cX\to \cY$ maps the closed point $x\in \cX(k)$ to a closed point $y\in \cY(k)$, then there is an induced map $\cN_x\to \cN_y$ on normal stacks. 

If $X=\Spec(A)$ is an affine scheme over $k$ and $G$ is a linearly reductive algebraic group, then $X/G\to \Spec(A^G)$ is a good moduli space. We use the notation $X\git G\coloneqq \Spec(A^G)$ for the affine GIT quotient.

If $\cX$ and $\cY$ are algebraic stacks admitting good moduli spaces $\cX\to X$ and $\cY\to Y$, then a morphism $f\colon \cX\to \cY$ is said to be \textit{strongly étale} if the induced morphism $X\to Y$ is étale and the square
\[\begin{tikzcd}
	\cX & \cY \\
	X & Y
	\arrow["f", from=1-1, to=1-2]
	\arrow[from=1-2, to=2-2]
	\arrow[from=2-1, to=2-2]
	\arrow[from=1-1, to=2-1]
\end{tikzcd}\]
is Cartesian (see \cite[Definition 3.13]{alper-hall-rydh-theetalelocalstructureofalgebraicstacks}).

\subsection{Slices for smooth morphisms of smooth stacks}

\begin{thm}\label{proposition: relative smooth etale slice}
Let $f\colon \cX\to \cY$ be a smooth $k$-morphism between quasi-separated smooth algebraic stacks over $k$ with affine stabilizers. Let $x\in \cX(k)$ be a closed point such that $y\coloneqq f(x)$ is closed. Suppose the stabilizers $G_x$ and $G_y$ are linearly reductive and that $G_y$ is smooth. Then there is a commutative diagram of pointed stacks
\begin{equation}
\begin{tikzcd}[ampersand replacement=\&] \label{diagram: local form all smooth stack}
	{(\cN_x,0)} \& {(\cU,u)} \& {(\cX,x)} \\
	{(\cN_y,0)} \& {(\cV,v)} \& {(\cY,y)}
	\arrow["f", from=1-3, to=2-3]
	\arrow[from=1-2, to=2-2]
	\arrow["h"', from=1-1, to=2-1]
	\arrow["{a_1}", from=1-2, to=1-3]
	\arrow["{b_1}"', from=2-2, to=2-3]
	\arrow["{b_2}", from=2-2, to=2-1]
	\arrow["{a_2}"', from=1-2, to=1-1],
\end{tikzcd}
\end{equation}
where 
\begin{enumerate}
    \item the maps $a_2, b_2$ are affine and strongly étale,
    \item the maps $a_1, b_1$ are étale and induce isomorphisms $G_u\stackrel{\sim}\to G_x$ and $G_v\stackrel{\sim}\to G_y$ of stabilizer groups, and
    \item the map $h$ is the one induced by  $f$  on normal spaces. 
\end{enumerate}  
In addition, if $\cX$ (resp. $\cY$) has affine diagonal, then $a_1$ (resp. $b_1$) can be arranged to be affine.
Moreover, if $\cX$ (resp. $\cY$) has affine diagonal and a good moduli space, then it can be arranged that $a_1$ (resp. $b_1$) is affine and strongly étale.
\end{thm}
\begin{remark}
    Note that $\cU$ and $\cV$ have good moduli spaces because $\cN_x$ and $\cN_y$ do and $a_2$ and $b_2$ are affine \cite[Lemma 4.14]{alper-good-moduli}.
\end{remark}
\begin{proof}
Since the result is local at the source and target,
     we may assume that $\cX$ and $\cY$ are quasi-compact. 
     By \cite[Thm. 1.2]{alper-hall-rydh-etale-slice}, we can find $(\cV,v)$, $b_1$ and $b_2$ as in
     \eqref{diagram: local form all smooth stack}, with $b_2$ affine and strongly étale. Forming the fiber product $\cZ\coloneqq \cX\times_\cY \cV$, we get a diagram
\[\begin{tikzcd}
	& {(\cZ,z)} & {(\cX,x)} \\
	{(\cN_y,0)} & {(\cV,v)} & {(\cY,y),}
	\arrow["f", from=1-3, to=2-3]
	\arrow["{f'}"', from=1-2, to=2-2]
	\arrow["{b_1}"', from=2-2, to=2-3]
	\arrow["{b_1'}", from=1-2, to=1-3]
	\arrow["\ulcorner"{anchor=center, pos=0.125}, draw=none, from=1-2, to=2-3]
	\arrow["{b_2}", from=2-2, to=2-1]
\end{tikzcd}\]
where $b_1'$ is étale and induces an isomorphism of stabilizer groups at $z$.
     Replacing $\cX$ by $\cZ$ and $f$ by the composition $\cZ\to \cV\to \cN_y$, we may assume that $(\cY,y)=(\cN_y,0)$. In order to construct diagram \eqref{diagram: local form all smooth stack}, we are reduced to producing a dashed étale roof $(\cU,u)$ fitting into the following commutative diagram
\[\begin{tikzcd}
	& {(\cU,u)} \\
	{(\cN_x,0)} && {(\cX,x)} \\
	& {(\cN_y,0)}
	\arrow["f", from=2-3, to=3-2]
	\arrow["h"', from=2-1, to=3-2]
	\arrow["{a_1}", dashed, from=1-2, to=2-3]
	\arrow["{a_2}"', dashed, from=1-2, to=2-1]
\end{tikzcd}\]
    with $a_2$ affine and strongly étale and with $a_1$ étale and inducing an isomorphism of stabilizer groups at $u$ and $x$.

    The point $x$ defines a closed immersion $BG_x\to \cX$ cut out by an ideal sheaf $\cI$. Let us denote by $\cX^{[n]}$ the $n$-th thickening of $x$, that is the closed substack of $\cX$ with ideal sheaf $\cI^{n+1}$. Similarly, we denote by $\cN_x^{[n]}$  and $\cN_y^{[n]}$ the $n$-th order thickenings of $0$ in $\cN_x$ and $\cN_y$, respectively. We begin by showing the following
    
    \begin{claim} There is an isomorphism $\cN_x^{[1]}\cong \cX^{[1]}$ fitting in a commutative diagram
\begin{equation}\label{diagram: square compatible isos thickenings} 
\begin{tikzcd}[ampersand replacement=\&,sep=scriptsize]
	{\cN_x^{[1]}} \&\& {\cX^{[1]}} \\
	{\cN_x} \&\& \cX \\
	\& {\cN_y}
	\arrow["\sim", from=1-1, to=1-3]
	\arrow[from=1-1, to=2-1]
	\arrow[from=1-3, to=2-3]
	\arrow["h"', from=2-1, to=3-2]
	\arrow["f", from=2-3, to=3-2]
\end{tikzcd}\end{equation}
\end{claim}
\noindent    \textit{Proof of the {\bf Claim}.} First note that we have canonical isomorphisms $BG_x\cong \cN_x^{[0]}$ and $BG_x\cong \cX^{[0]}$ over $\cN_y$.
Now consider the commutative diagram of solid arrows
\begin{equation}\label{diagram: first thickening}\begin{tikzcd}[ampersand replacement=\&,sep=scriptsize]
	\& {\cX^{[0]}} \\
	{\cX^{[1]}} \&\& {BG_x} \\
	\& {\cN^{[0]}_y} \\
	 {\cN_y^{[1]}} \&\& {BG_y,}
	\arrow[from=1-2, to=3-2]
	\arrow["\sim", no head, from=1-2, to=2-3]
	\arrow["c", from=2-3, to=4-3]
	\arrow[from=2-1, to=4-1]
	\arrow[from=1-2, to=2-1]
	\arrow[from=3-2, to=4-1]
	\arrow["\sim", no head, from=3-2, to=4-3]
	\arrow[from=4-1, to=4-3]
	\arrow[dashed, from=2-1, to=2-3]
\end{tikzcd}\end{equation}
where all vertical arrows are the ones induced by $f$. We shall produce a dashed arrow making the diagram commutative. By \cite[Thm. 8.5]{pridham-higher-stacks-presentation}, the obstruction to the existence of such a dashed arrow lives in $\Ext^1_{BG_x}(\bL_c,\cI/\cI^2)$. Since $G_y$ is smooth, the morphism $c$ is smooth, so the cotangent complex $\bL_c$ is concentrated in degrees $[0,1]$. Moreover, since $G_y$ is linearly reductive, the category $\QCoh(BG_y)$ of $G_y$-representations is semisimple. It follows that
\[\Ext^1_{BG_x}(\bL_c,\cI/\cI^2)=\Hom_{BG_x}(\bL_c,\cI/\cI^2[1])=0,\]
and a desired dashed arrow as in (\ref{diagram: first thickening}) exists. 

Since $\cX^{[0]}\to \cX^{[1]}$ is a surjective closed immersion and $\cX^{[0]}\to BG_x$ is representable, we have that $\cX^{[1]}\to BG_x$ is representable too. Since $\cX^{[0]}\to \cX^{[1]}$ is a square-zero extension and $\cX^{[0]}\to BG_x$ is affine, by dévissage the map $\cX^{[1]}\to BG_x$ is cohomologically affine. By Serre's criterion for algebraic spaces \cite[\href{https://stacks.math.columbia.edu/tag/07V6}{Tag 07V6}]{stacks-project}, we have that $\cX^{[1]}\to BG_x$ is affine. Diagram (\ref{diagram: first thickening}) induces compatible splittings $d_1$ and $d_2$ in the diagram
\[\begin{tikzcd}[ampersand replacement=\&]
	0 \& {N_x^\vee} \& {\cO_{\cX^{[1]}}} \& {\cO_{BG_x}} \& 0 \\
	0 \& {c^*N_y^\vee} \& {c^*\cO_{\cN_y^{[1]}}} \& {c^*\cO_{BG_y}} \& 0
	\arrow[from=2-1, to=2-2]
	\arrow[from=2-2, to=2-3]
	\arrow[from=2-3, to=2-4]
	\arrow[from=2-4, to=2-5]
	\arrow[from=1-4, to=1-5]
	\arrow[from=1-3, to=1-4]
	\arrow[from=1-2, to=1-3]
	\arrow[from=1-1, to=1-2]
	\arrow[from=2-2, to=1-2]
	\arrow[from=2-3, to=1-3]
	\arrow[from=2-4, to=1-4]
	\arrow["{d_1}"', curve={height=12pt}, from=1-4, to=1-3]
	\arrow["{d_2}"', curve={height=12pt}, from=2-4, to=2-3]
\end{tikzcd}\]
of short exact sequences in $\QCoh(BG_x)$. This gives an isomorphism $\cX^{[1]}\cong \cN^{[1]}_x=\Spec_{BG_x}(\cO_{BG_x}\oplus \varepsilon N_x^\vee)$, where $\varepsilon^2=0$, commuting with the maps into $\cN^{[1]}_y$.
The {\bf Claim} is thus established.

By \cite[Prop 7.18(2)]{alper-hall-rydh-theetalelocalstructureofalgebraicstacks}, there is an affine and strongly étale neighborhood 
\[a_2\colon (\cU,u)\to (\cN_x,0)\] 
and a morphism $a_1\colon (\cU,u)\to (\cX,x)$ over $\cN_y$ extending the map $\cN_x^{[1]}\xrightarrow{\sim}\cX^{[1]}\to \cX$ from the \textbf{Claim}.

Arguing as in the proof of \cite[Thm. 1.2]{alper-hall-rydh-etale-slice}, by smoothness of $\cX$ at $x$ and since $a_1$ induces an isomorphism $\cU^{[1]}\cong \cX^{[1]}$ of first order thickenings at $u$ and $x$, it follows from \cite[Prop. A.8]{alper-hall-rydh-etale-slice} that $a_1$ also induces isomorphisms $\cU^{[n]}\to \cX^{[n]}$ of $n$th order thickenings at $u$ and $x$ for all $n$. Therefore $a_1$ is étale at $u$ and thus there is an open neighborhood $\cU_1$ of $u$ in $\cU$ such that $a_1\vert_{\cU_1}$ is étale. Since $u$ is a closed point of $\cU$, there is another open neighborhood $\cU_2$ of $u$ that is saturated with respect to the good moduli space of $\cU$ and is contained in $\cU_1$. Therefore, after shrinking $\cU$ we may assume that $a_1$ is étale, while keeping $a_2$ affine and strongly étale. The map $a_1$ induces isomorphism of stabilizer groups at $u$ because $\cU^{[0]}\to \cX^{[0]}$ is an isomorphism of residual gerbes. We have proved the first part of the Theorem.

If $\cX$ has affine diagonal, then by \cite[Prop. 3.2]{alper-hall-rydh-etale-slice} there is an open substack  
$\cU'\subset \cU$ containing $u$, saturated with respect to the good moduli space $\pi\colon \cU\to U$ of $\cU$, such that the composition $\cU'\to \cU\to \cX$ is affine. Since $\pi^{-1}\pi(\cU')=\cU'$, the map $\cU'\to \cU$ is in particular strongly étale. If moreover $\cX$ has a good moduli space, then by Luna's fundamental lemma for stacks \cite[Thm. 3.14]{alper-hall-rydh-theetalelocalstructureofalgebraicstacks} there is a $\pi$-saturated open neighborhood $\cU''$ of $u$ in $\cU'$ such that $\cU''\to \cX$ is strongly étale. The same arguments work for $\cY$.\end{proof}

\begin{remark}
We thank David Rydh for pointing out that \cite[Prop 7.18(2)]{alper-hall-rydh-theetalelocalstructureofalgebraicstacks} can be used to shorten our original proof of \Cref{proposition: relative smooth etale slice}. The original argument used Tannaka duality and Artin approximation, modifying the arguments in \cite{alper-hall-rydh-etale-slice} to ensure that the relevant morphisms commute with the maps to the base algebraic stack $\cN_y$. We refer to the first preprint version of this article \cite[Thm. 2.1]{decataldo2023relative} for details.
\end{remark}

In the setting of \Cref{proposition: relative smooth etale slice}, suppose that the base $\mathcal{Y}=S$ is an algebraic space. Let $N_{f,x}$ denote the kernel of $N_x \to N_{f(x)}$. By the linear reductivity of $G_x$, we have an isomorphism of $G_x$-representations $N_{x} \cong N_{f,x} \oplus T_{S, f(x)}$, where $G_x$ acts trivially on the tangent space $T_{S, f(x)}$. 

\begin{coroll}[Equisingularity of good moduli spaces of smooth stacks] \label{cor: equisingularity}
Let $f: \cX \to S$ be a smooth algebraic stack with affine diagonal over a smooth $k$-algebraic space $S$. Suppose that $\cX$ admits a good moduli space  $\overline{f}: M \to S$. Let $x \in \cX(k)$ be a closed point with corresponding point $x \in M(k)$. Then there is a commutative diagram of pointed algebraic spaces
\begin{equation} 
\begin{tikzcd}[ampersand replacement=\&]\label{diagram: local form all smooth}
	(T_{S,\overline{f}(x)}\times (N_{f,x}\git G_x), \, (0,0)) \ar[d, "h"] \& \ar[l, "a_2", labels=above] (U,u) \ar[d, "g"] \ar[r, "a_1"] \& (M,x) \ar[d, "\overline{f}"] \\
	(T_{S, \overline{f}(x)}, 0) \& \ar[l, "b_2", labels=above] (V,v) \ar[r, "b_1"] \& (S, \overline{f}(x))
\end{tikzcd}
\end{equation}
where $a_1, a_2$ and $b_1, b_2$ are \'etale.
\end{coroll}

\begin{proof}
    The hypothesis that $\cX$ admits a good moduli space implies that the closed $k$-point $x \in \cX(k)$ has linearly reductive stabilizer \cite[Prop. 12.4]{alper-good-moduli}. Hence, the corollary follows by taking good moduli spaces in  diagram \eqref{diagram: local form all smooth stack} of \Cref{proposition: relative smooth etale slice},  where the horizontal morphisms can be arranged to be strongly \'etale because $f$ has affine diagonal. 
\end{proof}

In plain words, \Cref{cor: equisingularity} says that the singularities of a family $\overline{f}: M \to S$ obtained as the good moduli space of a smooth $S$-stack are 
\'etale locally constant.

\begin{remark}
    David Rydh has explained to us in a private communication that it is also possible to deduce \Cref{cor: equisingularity} from \cite[Cor. 8.7(8)]{alper-hall-rydh-theetalelocalstructureofalgebraicstacks}.
\end{remark}

\section{Cohomology in families} \label{section: cohomology in families}

In this section, we shall prove results about the (intersection) cohomologies of fibers of a locally isotrivial morphism in the following sense.
\begin{defn} \label{defn: locally isotrivial}
    Let $f: \cM \to \cS$ be a morphism of finite type algebraic stacks over $k$. We say that $f$ is locally isotrivial if there is a smooth surjective morphism $\cY \to \cS$ such that the base-change $f_{\cY}: \cM_{\cY} \to \cY$ fits into a commutative diagram
    \begin{equation} \label{diagram: local isotriviality}
\begin{tikzcd}[ampersand replacement=\&]
	{\cT \times \cF} \& {\cU} \& {\cM_{\cY}} \\
	{\cT} \& {\cV} \& {\cY}
	\arrow["f_{\cY}", from=1-3, to=2-3]
	\arrow[from=1-2, to=2-2]
	\arrow["p_{\cT}"', from=1-1, to=2-1]
	\arrow["{a_1}", from=1-2, to=1-3]
	\arrow["{b_1}"', from=2-2, to=2-3]
	\arrow["{b_2}", from=2-2, to=2-1]
	\arrow["{a_2}"', from=1-2, to=1-1],
\end{tikzcd}
\end{equation}
where $p_{\cT}$ is the first projection, $a_1,b_1$ are surjective \'etale morphisms, and $a_2, b_2$ are \'etale.
\end{defn}

\Cref{cor: equisingularity} says that good moduli spaces of smooth stacks over a base are locally isotrivial.

\subsection{Preliminaries on \texorpdfstring{$\ell$}{l}-adic sheaves} \label{preliminaries l-adic}
For the next subsections, we fix a prime 
number $\ell$
distinct from the characteristic of the algebraically closed ground field $k,$ and we work with $\overline{\mathbb{Q}}_{\ell}$-sheaves for the \'etale topology.

Let $\cT$ be an algebraic stack which is quasi-separated and of finite type over an algebraically closed field,
or over the spectrum of a strictly Henselian discrete valuation ring. Let $D^b_c(\cT, \overline{\mathbb Q}_\ell)$ be the $\overline{\mathbb Q}_\ell$-constructible derived category endowed with the usual six functors $(f^*, Rf_*, Rf_!, f^!, \otimes^{L}, R{\cH}om)$. Standard references for the construction
of the ``$\ell$-adic constructible derived" category and of the six functors are \cite{ekedahl-adic,proetale} in the context of schemes, and \cite{laszlo-olsson-adic-sheaves} in the context of quasi-separated algebraic stacks.

Let $\cM$ be a quasi-separated finite type algebraic stack over the algebraically closed field
$k.$ We use the notion of perverse $t$-structure defined in \cite{laszlo-olsson-perverse-sheaves}. Let $IC_{\cM}$ be the perverse intersection cohomology complex of $\cM$ with coefficients in $\overline{\mathbb{Q}}_{\ell}$, defined as the direct sum of the perverse intersection cohomology complexes of the irreducible components of $\cM$.
If $\cM$ is irreducible, then we define the topologist's intersection complex $\mathcal{IC}_{\cM}$ 
of $\cM$ by setting
$\mathcal{IC}_{\cM} :=  IC_{\cM}[-\dim(\cM)].$
If $\cM$ is
not irreducible with irreducible components $\cM_i$, then we define the topologist's intersection complex of $\cM$ by setting $\mathcal{IC}_{\cM}:= \oplus_i \mathcal{IC}_{\cM_i}.$

\begin{notn}
We denote by $H^{\bullet}(\cM, \overline{\mathbb{Q}}_{\ell})$ the cohomology ring of $\cM$. Similarly, we write $I\!H^{\bullet}(\cM, \overline{\mathbb{Q}}_{\ell})$ for the intersection cohomology of $\cM$. 
\end{notn}

We will also make use of the theory of nearby and vanishing cycle functors, denoted $\psi$ and $\phi$ respectively, for morphisms of schemes. 
Standard references for nearby and vanishing cycle functors, $\psi$ and $\phi$, are \cite[Exp. XIII]{SGA7-2} and \cite{illusie-monodromie, illusie-perversite-et-variation}. See also 
\cite[Appendix \S5.2]{decataldo-zhang-completion}, and \cite[\S2.1]{decataldo-cambridge}
for a partial list of the properties concerning 
nearby and vanishing cycle functors. The facts about nearby and vanishing cycle functors we use freely in this paper are:
\begin{enumerate}[(1)]
    \item the existence of the functorial distinguished 
triangles of functors $i^* [-1] \to \psi [-1] \to \phi  \rightsquigarrow,$ and 
$\phi \to \psi [-1] \to i^! [1] \rightsquigarrow$;
for the first one, see \cite[Exp. XIII, 2.1.2]{SGA7-2}; for the second one, apply Verdier Duality to the first one: see \cite[Thm. 4.2]{illusie-monodromie} for the compatibility of nearby cycles with duality; see
\cite[Cor. 0.2]{phi-duality}
for the compatibility of vanishing cycles with duality;

\item the identities $Rf_*\psi= \psi Rf_*$ and $Rf_* \phi = \phi Rf_*$ for a proper morphism $f$
(cf. \cite[Exp. XIII, 2.1.7]{SGA7-2});

\item the identities $f^*\psi= \psi f^*$ and $f^* \phi = \phi f^*$ for a smooth morphism $f$
(cf. \cite[Exp. XIII, 2.1.7]{SGA7-2});

\item  the compatibility of
$\psi$ with cup products:
$\psi (F) \otimes^L \psi (G) \stackrel{\sim}\to \psi (F\otimes^L G)$
(cf. \cite[Thm. 4.7]{illusie-monodromie}).
\end{enumerate}

\begin{remark}\label{phi and ag sp}
While the six-functor formalism is available for algebraic spaces and stacks, we could not locate a suitable reference for the nearby and vanishing cycle functors. In view of this, we do not state \Cref{supp coeff} in the generality of proper DM morphisms, and instead we restrict ourselves to proper schematic morphisms.
\end{remark}

 \subsection{\'Etale cohomology of locally isotrivial families}\label{ac field}
 In this subsection, we work with $\overline{\mathbb{Q}}_{\ell}$-sheaves. Let $f:\cM \to \cS$ be a proper morphism between algebraic stacks which are of finite type and with affine diagonal over the algebraically closed field
$k.$

By the Decomposition Theorem \cite{bbdg, shenghao-decomposition-stacks} there is a decomposition
\begin{equation}\label{DT}
Rf_* (\mathcal{IC}_{\cM}) \cong \bigoplus_{i} {}^p\mathcal{H}^i(Rf_* (\mathcal{IC}_{\cM}))[-i]
\end{equation}
of the direct image into the finite direct sum of its shifted perverse cohomology sheaves.
Up to shift, the non-trivial summands are semisimple perverse sheaves on $\cS$ whose 
simple components are intersection complexes
$\mathcal{IC}_{\cZ}(L)[dim(\cZ)]$ of closed
irreducible substacks $\cZ \subset \cS$ (the supports) with coefficients in irreducible lisse sheaves $L$ (the simple coefficients) on a suitable
open dense substack $\cZ^o\subseteq \cZ.$

In general, the determination of the supports and the coefficients  is a difficult task.
The following is one of the main theorems of this subsection, which provides some conditions that ensure that the only support is the base $\cS$ and the coefficients
are the lisse sheaves associated with the (intersection) cohomology
groups of the fibers.
\begin{thm}[Supports and coefficients]\label{supp coeff}
Let $f: \cM \to \cS$ be a proper schematic morphism between algebraic stacks which are of finite type and with affine diagonal over $k$. Suppose that $f$ is locally isotrivial in the sense of \Cref{defn: locally isotrivial}, and that $\cS$ is smooth and connected. Then the following hold:
\begin{enumerate}
\item[(a)]
For every $i,$ the direct image sheaves
$R^i f_* (\overline{\mathbb{Q}}_{\ell})_{\cM}$ and
$R^i f_* \mathcal{IC}_{\cM}$
are lisse on $\cS$. 

\item[(b)]
We have the identity of perverse sheaves on $\cS$
\[
(R^if_* (\mathcal{IC}_{\cM}))[{\rm dim} (\cS)] =
{}^p\mathcal{H}^{\dim (\cS)+i} (Rf_* (\mathcal{IC}_{\cM})).
\]
\item[(c)]
We have an isomorphism
in $D^b_c(\cS, \overline{\mathbb Q}_\ell)$
\[
Rf_* \mathcal{IC}_{\cM} \simeq \bigoplus_{i\geq 0}
R^if_* \mathcal{IC}_{\cM}[-i].
\]
\end{enumerate}
\end{thm}

Our next goal in this subsection is to prove \Cref{supp coeff} and the upcoming \Cref{iso coh fibers}. We need the following seemingly standard three lemmata, for which we could not find suitable references.
 
\begin{lemma}\label{ic field}
 The formation of the intersection complex $\mathcal{IC}_{X_E}$ of a scheme of finite type $X_E$ over a field $E$ is compatible with field extensions $L/E,$
 i.e. we have a natural isomorphism
 $\mathcal{IC}_{X_L}= u^*
 \mathcal{IC}_{X_E},$ where $u:X_L:= X_E\times_E L 
 \to X_E$ is the natural projection.
 \end{lemma}

 \begin{proof}
 This follows from  the definition of intersection complexes as intermediate extensions (cf. \cite[Déf. 1.4.22]{bbdg}) and  from the standard fact that the functor $u^*:D^b_c(X_E, \overline{\mathbb Q}_\ell) \to D^b_c(X_L, \overline{\mathbb Q}_\ell)$
 is $t$-exact for the respective middle-perversity $t$-structures.
 
In order to establish this $t$-exactness, one can argue as follows.
Recall that the perverse $t$-structure $({}^p D^b_c(X_E)^{\leq 0}, {}^p D^b_c(X_E)^{\geq 0})$ on $X_E$ is obtained by gluing the perverse 
$t$-structures on any partition $X_E = U_E \cup F_E$
into complementary open and closed subsets
(cf. \cite[Thm. 1.4.10]{bbdg}). We need to show that
the functor $u^*$ preserves the subcategories of type
$({}^p D^b_c(X_?)^{\leq 0}, {}^p D^b_c(X_?)^{\geq 0})$. By Noetherian induction, we can assume that the desired $t$-exactness of $u^*$ holds when $X_E$ is replaced by $U_E$ and by $Z_E$, respectively. Consider the Cartesian diagram in \cite[pg.155]{bbdg}  (with $U=S=Spec(E)$ in loc.cit) and its  Cons\'equence in loc.cit., i.e.  the compatibility of $u^*$ with the formation of the six functors (we only need $j^*$ for $j$ an open immersion, and $i^*$ and $i^!$ for  $i$ a closed immersion).
This compatibility yields the desired preservation
under $u^*$
of the subcategories of type
$({}^p D^b_c(X_?)^{\leq 0}, {}^p D^b_c(X_?)^{\geq 0})$.
\end{proof}

\begin{lemma} \label{lemma: fiber of IC vs IC of fiber}
Let $f: M \to S$ be a morphism of finite type schemes over $k$. Suppose that $f$ is locally isotrivial in the sense of \Cref{defn: locally isotrivial} and that $S$ is smooth.
For every geometric point $s \to S$, let $\iota_s: M_s \to M$ denote the corresponding geometric fiber. Then, we have a canonical isomorphism $\iota_s^* \mathcal{IC}_{M} = \mathcal{IC}_{M_s}$. 
\end{lemma}
\begin{proof}
  Since the intersection complex $\mathcal{IC}_{M}$ is compatible with extension of the ground field (\Cref{ic field}), we can replace $k$ with the residue field of $s$ and assume without loss of generality that $s$ is a closed $k$-point in $S.$ Since $f$ is locally isotrivial, it is in particular flat. Therefore every irreducible component of $M$ restricts  to a union of irreducible components of $M_s$. By the definition of $\mathcal{IC}_M$, we may restrict to each irreducible component of $M$ and assume without loss of generality that $M$ is integral.

    Consider the local model commutative diagram
    \eqref{diagram: local isotriviality}, where we may assume that everything is a scheme, and so we shall use roman letters $V,U, F, T$ for the schemes appearing in the diagram. Choose $k$-points $0 \in T(k)$ and $v \in V(k)$ such that $b_1(v) =s$ and $b_2(v) =0$.
    We take the fibers over the points
    $0,$ $v$ and $s$ and obtain the commutative diagram (with some decorations omitted)
    \begin{equation} 
    \begin{tikzcd}[ampersand replacement=\&]\label{diagram fiber}
	 F  \ar[d, "\iota_0"] \& \ar[l, "a'_2", labels=above] 
 U_v \ar[d, "\iota_v"] \ar[r, "a'_1"] \& M_s \ar[d, "\iota_s"] \\
    T \times F \& \ar[l, "a_2", labels=above] U \ar[r,"a_1"] \& M
\end{tikzcd}
\end{equation}
where $a_1, a_2$ and $a'_1, a'_2$ are \'etale, and $a_1$ and $a_1'$ are surjective.
By working with the commutative square on the left-hand side
of (\ref{diagram fiber}), 
we obtain the following chain of natural isomorphisms
\[ \iota_v^* \mathcal{IC}_U=
\iota_v^* a_2^* \mathcal{IC}_{T \times F} =(a_2')^*\iota_0^*\mathcal{IC}_{T \times F}= (a'_2)^*\mathcal{IC}_{F} = \mathcal{IC}_{U_v},
\]
where: the first natural identification is because $a_2$ is \'etale;
the second is because the diagram is commutative;
the third is because we can factor the identity as  ${\rm id}_{F}: F \xrightarrow{\iota_0}T \times F \xrightarrow{p_F} F,$ where the projection morphism satisfies $p_F^*\mathcal{IC}_{F} = \mathcal{IC}_{T\times F}$ because it is smooth; the fourth is because $a_2'$ is \'etale.

By working   with the commutative square on the right-hand side
of (\ref{diagram fiber}) the same way we have done on the left-hand side, we see that
\[
(a_1')^*\iota_s^* \mathcal{IC}_M  = 
\iota_v^* a_1^* \mathcal{IC}_M =  \iota_v^* \mathcal{IC}_U
= \mathcal{IC}_{U_v}.
\]
Since $a_1'$ is \'etale, we also have that $(a_1')^* \mathcal{IC}_{M_s}= \mathcal{IC}_{U_v}.$

By our reduction, we are assuming that $M$ is integral, and therefore every connected component of $T \times F$ is integral. This in particular implies that every connected component of $F$ is integral, and therefore the same holds for $M_s$.
By considering the perverse versions of the intersection complexes,
we have a natural isomorphism of perverse shaves
\[
(a_1')^* IC_{M_s}= IC_{U_v} 
= (a_1')^*\iota_s^* IC_M [\dim{S}-\dim{M}].
\]
Since being a perverse intersection complex is a local condition in the \'etale topology, and the morphism $a_1'$ is \'etale and surjective, we deduce that 
$\iota_s^* IC_M [\dim{S}-\dim{M}]$ is an intersection complex
on $M_s$
with possibly twisted coefficients on $M_s.$ By restricting to the smooth locus of $M$ where the intersection complex is the constant sheaf,
we see that the identification sends the unit section to the unit section, so that the twisted coefficients in question are in fact constant of rank one. 
The conclusion  of the lemma follows.
\end{proof}

\begin{remark}
    For \Cref{IC and product} below, we use O. Gabber's rectified middle perversity structures for schemes of finite type
over a DVR $T;$ see \cite{illusie-monodromie,illusie-perversite-et-variation}; for a summary and complements, see \cite[\S5.2]{decataldo-zhang-completion}. In short, one glues
the middle-perversity $t$-structure on the closed fiber with the   middle-perversity $t$-structure on the generic fiber shifted by $1.$ Just as in the case of schemes of finite type over  a field, one has, in the context  of schemes of finite type over a DVR, the notion intermediate extension of a perverse sheaf from an open set to a bigger open set, and thus the usual notion of intersection complex on such schemes.
\end{remark}

The following lemma, for which we could not locate a reference,  is standard (K\"unneth Formula for intersection cohomology) when $T$ is a scheme of  finite type over the field $k.$ 

\begin{lemma}\label{IC and product}
Let $N$ be a scheme of finite type over an algebraically closed field $k.$
Let $T$ be the spectrum of a strictly Henselian DVR over $k$ and with residue field $k.$ Then
$\mathcal{IC}_{N\times T} = p_N^*
\mathcal{IC}_N,$ where $p_N: N\times T \to N$
is the projection.
\end{lemma}
\begin{proof}
Without loss of generality, we assume that $N$ is
integral. We use the perverse version
$IC$ of
intersection cohomology.
Let $i: s\to  T \leftarrow  \eta :j$ be the closed and open point, respectively. 
Let $p_T: N\times T \to T$ be the projection and let $N=N_s$ and  $N_\eta$ be the special and generic fiber respectively. We have the following commutative diagram with Cartesian squares
\begin{equation} 
    \begin{tikzcd}[ampersand replacement=\&]\label{diagram ic product}
    \& N \&
    \\
	N \ar[ru, "id"]  \ar[r, "i"] \ar[d]  \&
 N\times T \ar[d, "p_T"] \ar[u, "p_N", labels=right] \ar[d]  \& N_\eta  \ar[d] \ar[l, "j" , labels=above] \ar[ul, "u", labels=above] \\
	s  \ar[r, "i"] \& T  \& \eta \ar[l, "j" , labels=above]
\end{tikzcd}
\end{equation}

 We have morphisms $i^*[-1] (p^*_N IC_N[1])
\to \psi[-1]  (p^*_N IC_N[1]) \to i^! [1] (p^*_N IC_N[1]).$ A cone for  either morphism is $
\phi (p_N^*IC_N[1]),$ which is zero in view of 
\cite[pg.243, Cor. 2.16]{SGA4.5}. Hence,
we have the natural identification 
$i^*[-1] (p^*_N IC_N[1]) = 
i^! [1] ( p^*_N IC_N[1]).$  

Since $i^* (p_N^* IC_N)  =i^! (p_N^* IC_N) [2] = IC_N$
and $j^* (p_N^* IC_N) = IC_{N_\eta}$
(cf. \Cref{ic field}),
one verifies (cf. \cite[(55)]{decataldo-zhang-completion})  that $p_N^*(IC_N[1])$ is perverse on $N\times T$ for the rectified $t$-structure.

We have a distinguished attaching triangle
\[
i_* IC_N [-1] =
i_* i^! p_N^* (IC_N [1]) \to p_N^* (IC_N [1]) \to Rj_* j^* p_N^*
(IC_N [1]) = Rj_*IC_{N_\eta}[1] \rightsquigarrow
\]
Note that all three complexes $p_N^*IC_N [1]$, $i_*(IC_N)$ and $Rj_*(IC_{N_{\eta}})[1]$ are perverse, the last two by the $t$-exactness of the functors $i_*, Rj_*[1]$ and the perversity of $IC_N, IC_{N_{\eta}}$. By taking the long exact sequence of perverse cohomology sheaves, we get the short exact sequence of perverse sheaves on $N \times T$
\[
0 \to p_N^* (IC_N [1]) \to Rj_* (IC_{N_{\eta}}[1])
\to i_* IC_N \to 0.
\]

Since $IC_{N\times T}$ is the smallest sub perverse sheaf
of $Rj_* (IC_{N_{\eta}}[1])$ extending 
$IC_{N_{\eta}}[1]$ (viewed as a rectified perverse sheaf on the open $N_\eta \subseteq N\times T$), we see that we have a natural short exact sequence of perverse sheaves on $N\times T$
\begin{equation}\label{poiqw}
0 \to IC_{N\times T} \to p_N^* (IC_{N}[1])
\to i_*Q_N \to 0,
\end{equation}
where $Q_N$ is a perverse sheaf on $N.$
Since $Q_N$ is supported on $N,$ we have that $\phi (i_* Q_N) = Q_N.$
We apply the exact functor $\phi$ 
to the short exact sequence (\ref{poiqw})
and, in view of the earlier-established  vanishing $\phi (p_N^*IC_N [1])=0,$ we deduce that 
$Q_N = \phi (i_* Q_N)=0.$
\end{proof}

\begin{fact}\label{trivial t str}
Let $S$ be a smooth irreducible scheme over 
the algebraically closed  field $k$, and let $K$
be  a complex with lisse cohomology sheaves.
Then we have a natural identification  $\mathcal{H}^i (K) [dim(S)] = {}^p\mathcal{H}^{dim(S)+i}(K).$
 In fact, by the
gluing construction
\cite[\S1.4]{bbdg}
of $t$-structures, on the dense stratum the $t$-structure is the standard one, up to dimensional shift.
\end{fact}

Now we are ready to conclude the proof of \Cref{supp coeff}.
\begin{proof}[Proof of \Cref{supp coeff}]
The conclusions of the theorem can be checked smooth locally on $\cS$. After base-changing to a smooth atlas $S \to \cS$, we are immediately reduced to the case when $f: M \to S$ is a proper morphism of schemes. In view of Fact \ref{trivial t str}, part (b)
follows from part (a). In view of the Decomposition Theorem   (\ref{DT}), part (c) follows from part (b).
We are left with proving part (a).

In what follows, we use implicitly the Proper Base Change Theorem for the morphism $f,$
so that the  formation of the direct image sheaves in question, simply denoted by $E$ in what follows,  are compatible with base change.
According to \cite[\href{https://stacks.math.columbia.edu/tag/0GKC}{Tag 0GKC}]{stacks-project},
we need to show that for every specialization 
$t \rightsquigarrow s$ of Zariski points in $S,$
and for any choice of geometric points
$\overline{s}$ over $s$ and $\overline{t}$ over $t,$
the resulting specialization morphism
$E_{\overline{s}} \to E_{\overline{t}}$ is an isomorphism. Without loss of generality, we can assume that $s \in S$ is a closed point and take $\overline{s}=s.$ 
By \cite[\href{https://stacks.math.columbia.edu/tag/0GJ6}{Tag 0GJ6}]{stacks-project}, we can furthermore replace $S$ with 
its strict localization $\widehat{S}$ at $s.$ 
Consider the  spectrum of a strictly Henselian DVR $\widetilde{S}$ mapping to $\widehat{S},$ 
and mapping its closed point $\widetilde{s}$ to $s$ and its generic point $\widetilde{t}$ to $t.$ 
We have the corresponding morphisms $\widetilde{S} \to \widehat{S} \to S$ which give rise  to 
 the corresponding morphisms obtained from $f$ via base change: $\widetilde{f}: \widetilde{M} \to \widetilde{S}$
and
$\widehat{f}: \widehat{M} \to \widehat{S}.$
Let $\overline{\widetilde{t}}$ be a geometric point over $\widetilde{t}.$  The specialization morphism 
$E_{\widetilde{s}} \to E_{\overline{\widetilde{t}}}$
coincides with the morphism $i^*E \to \psi(E)$
associated with $\widetilde{S}.$ We are thus reduced to showing  the triviality of the following vanishing cycle complexes 
\[\phi 
\left(R\widetilde{f}_* \left(\overline{\mathbb{Q}}_{\ell}\right)_{\widetilde{M}}\right) = 
\phi \left(R \widetilde{f}_* \mathcal{IC}_{\widetilde{M}}\right)  =0.
\]

We have $\phi R\widetilde{f}_* = R (\widetilde{f}_s)_*\phi$ because $\widetilde{f}$ is proper. Hence, we are reduced to proving that
\begin{equation}\label{to prove}
\phi 
\left(\left(\overline{\mathbb{Q}}_{\ell}\right)_{\widetilde{M}}\right) = 
\phi \left(\mathcal{IC}_{\widetilde{M}}\right)  =0.
\end{equation}

We have the local model expressed by diagram  \eqref{diagram: local isotriviality} for $f:M\to S$, where we may assume that all of the objects $T,U,V,F$ are schemes. The morphisms $b_1$ and $b_2$ in said diagram induce isomorphisms of strict localizations
at chosen corresponding points
$\widehat{T} \stackrel{\widehat{b_2} \simeq }\longleftarrow \widehat{V} \stackrel{\widehat{b_1} \simeq }\longrightarrow  \widehat{S}.$ It follows that to give a morphism $\widetilde{S} \to \widehat{S}$ is the same thing as giving a morphism
to either $\widehat{T}$ or $\widehat{V}.$ We base change the local model diagram first via the morphism $\widehat{S} \to S$ and then 
via $\widetilde{S} \to \widehat{S}$ and obtain the following commutative diagram:

\begin{equation} 
\begin{tikzcd}[ampersand replacement=\&]\label{bs cg}
	\widetilde{T} \times F \ar[d, "p_{\widetilde{T}}"] \& \ar[l, "\widetilde{a_2}", labels=above] \widetilde{U} \ar[d, "g"] \ar[r, "\widetilde{a_1}"] \& \widetilde{M} \ar[d, "\widetilde{f}"] \\
	\widetilde{T} \& 
 \ar[l, "\simeq", "\widetilde{b_2}"' ] 
 \widetilde{V} \ar[r, "\widetilde{b_1}", "\simeq"' ] \& \widetilde{S},
\end{tikzcd}
\end{equation}
where $\widetilde{a_1}$ and $\widetilde{a_2}$ are \'etale and $\widetilde{b_1}$ and $\widetilde{b_2}$ are isomorphisms.

In view of  Lemma \ref{IC and product},
the intersection complex $\mathcal{IC}_{F}$
pulls back via the  projection morphism $p_F$ to
the intersection complex $\mathcal{IC}_{\widetilde{T}\times F},$
which, in turn, pulls back via the \'etale morphism $\widetilde{a_2}$ to the intersection complex $\mathcal{IC}_{\widetilde{U}}.$
The intersection complex $\mathcal{IC}_{\widetilde{M}}$ also pulls back to  $\mathcal{IC}_{\widetilde{U}}.$

In order to simplify the notation,
for a scheme $Z$ we let  $K_Z$ denote either the intersection complex
${\mathcal IC}_Z,$ or the constant sheaf 
$\left(\overline{\mathbb{Q}}_{\ell}\right)_{Z}.$
We then have:
\[\widetilde{a_2}^* \left(K_{\widetilde{T} \times F}\right) = K_{\widetilde{U}} = \widetilde{a_1}^*\left(K_{\widetilde{M}}\right)\]

Since $K_{\widetilde{T} \times F}=p_T^*K_{F},$
by virtue of \cite[p.
243, Cor. 2.16]{SGA4.5} for the morphism $F \to \Spec(k)$, we have 
that $\phi_{p_{\widetilde{T}}}\left(K_{\widetilde{T} \times F}\right)=0.$

We use freely that the formation of vanishing cycles
``commutes" with smooth pull-back \cite[Exp. XIII, (2.1.7.2)]{SGA7-2}. We apply this repeatedly to the \'etale morphisms of type $\widetilde{a}.$ It follows that 
\[0= \widetilde{a_2}^* \phi_{p_{\widetilde{T}}}\left(K_{\widetilde{T} \times F}\right)=
\phi_{p_{\widetilde{T}}\widetilde{a_2}}\left(\widetilde{a_2}^* K_{ \widetilde{T} \times F}\right) = \phi_{\widetilde{b_2} g}\left(\widetilde{a_2}^* K_{ \widetilde{T} \times F}\right) = 
\phi_{\widetilde{b_2} g}\left(K_{\widetilde{U}}\right),\]
 so that, since $\widetilde{b_2}$
is an isomorphism, we have 
$\phi_{g}\left(K_{\widetilde{U}}\right)=0.$

We have 
\[\widetilde{a_1}^* \left(\phi_{\widetilde{f}}  \left(K_{\widetilde{M}}\right)\right) =
\phi_{\widetilde{f} \widetilde{a_1}} \left(\widetilde{a_1}^*K_{\widetilde{M}}\right) =\phi_{\widetilde{b_1} g} \left(\widetilde{a_1}^*K_{\widetilde{M}}\right) =
\phi_{\widetilde{b_1} g} \left(K_{\widetilde{U}}\right) =0,\]
where: the first and second  equalities hold  because $\widetilde{a_1}$ is \'etale; the third equality holds because $\widetilde{b_1}$ is an isomorphism
and we already know that 
$\phi_{g}\left(K_{\widetilde{U}}\right)=0.$
We have thus shown that the desired triviality
of vanishing cycles complexes (\ref{to prove})
holds locally on $\widetilde{M},$ so that it holds on $\widetilde{M}.$
This concludes the proof of part (a).
\end{proof}

\begin{coroll}[(Intersection) cohomology and fibers]\label{iso coh fibers}
Let $f\colon \cM \to \cS$ be a proper schematic morphism between quasi-separated finite type algebraic stacks over $k$. Suppose that $f$ is locally isotrivial (\Cref{defn: locally isotrivial}), and that $\cS$ is connected. Then the cohomology rings of any two geometric fibers  of $f$
are non-canonically isomorphic. Their intersection
cohomology groups, endowed with their natural  module over the respective cohomology rings structures, and with their natural intersection pairings, are also non-canonically isomorphic.
\end{coroll}
\begin{proof}
After base-changing to connected components of the atlas $S \to \cS$, we are immediately reduced to the case when $f\colon M \to S$ is a proper morphism of schemes.
We begin by showing the statement for the fibers over two closed points $s_1, s_2 \in S$.
By connectedness of $S$, there exists a connected reduced, possibly reducible algebraic curve on $S$ connecting the two points (to see this we may reduce to the case when $S$ is an affine and integral, where \cite[Lemma, pg.56]{mumford-abelian-varieties} applies). By taking the normalization of such a curve, we reduce to the case when $S$ is a smooth and connected curve. Ignoring the ring, module and pairing structures, the additive part of the corollary now follows from Theorem \ref{supp coeff}, proper base-change for $f$ and \Cref{lemma: fiber of IC vs IC of fiber}.

As to the ring, module  and pairing structures, 
it is enough to prove the desired conclusion
by replacing $s_2$ with the generic point $\eta$ of
the smooth curve, and replacing the curve by its strict localization
$\widetilde{S}$ at $s:=s_1$.  In this case, we need to prove that
the morphism induced in cohomology by $i^* \to \psi$ is compatible with the ring and module structures.

The nearby cycle functor $\psi$ preserves cup
products (see e.g. \cite[\S4.3]{illusie-monodromie}). Since $f$ is proper, upon taking stalks
(see \cite[Exp. XIII, (2.1.8.5)]{SGA7-2}) the morphism $i^* \to \psi$ induces
the specialization morphism on stalks as defined in \cite[\href{https://stacks.math.columbia.edu/tag/0GJ2}{Tag 0GJ2}]{stacks-project}. By the description of the
specialization morphism in terms of pulling back sections  as in \cite[\href{https://stacks.math.columbia.edu/tag/0GJ3}{Tag 0GJ3}]{stacks-project}, we see that the morphism $i^* \to \psi$ preserves cup products.

As to the module structure, we note that it can be described as follows. Let $\widetilde{S}$ be the strict localization of $S$ at $s$, with a choice of geometric generic point denoted by $\overline{\eta}$.
We then have that $i^* \to \psi$ induces isomorphisms
$H^*(M_s) = H^* (\widetilde{M}) \simeq H^*(M_{\overline{\eta}}),$ and similarly for intersection cohomology. Let $\alpha \in H^a (\widetilde{M})$
and let $\beta \in I\!H^b(\widetilde{M}).$
We view $\alpha$ as a morphism $\alpha': (\overline{\mathbb{Q}}_{\ell})_{\widetilde{M}}
\to (\overline{\mathbb{Q}}_{\ell})_{\widetilde{M}} [a];$
this morphism induces a morphism
$\alpha'': \mathcal{IC}_{\widetilde M} \to \mathcal{IC}_{\widetilde{M}} [a]$
such that the cup product $\alpha \cdot \beta = H^b(\alpha'')(\beta)$ is obtained by taking cohomology.
By plugging this construction into the morphism 
of functors $i^* \to  \psi,$ we see that
 the product $\alpha_s \cdot \beta_s$  is sent to the product  
$\alpha_{\overline{\eta}} \cdot \beta_{\overline{\eta}}.$

Finally, we consider the pairings on intersection cohomology.  We place ourselves
again over $\widetilde{S}.$ 
The pairing on  intersection cohomology
arises from the map in cohomology induced by the self-duality isomorphism $\iota: IC \xrightarrow{\sim} {\mathbb D} IC$ (here, $\mathbb D$ denotes the Verdier Duality Functor and, in order to avoid dimensional shifts, we consider  perverse intersection cohomology complexes $IC$, which are self-dual).
We start with $\iota_M: IC_M \to {\mathbb D}_M IC_M$ and apply $i^*\to \psi.$ In order to conclude, it remains to observe that: 
\begin{enumerate}[(1)]
    \item $i^*IC_M  = IC_{M_s}[1]$ and $\psi (IC_M)  = IC_{M_{\overline{\eta}}}[1];$
    \item  $i^* {\mathbb D}_M IC_M =
{\mathbb D}_{M_s} i^! IC_M = {\mathbb D}_{M_s} 
(i^* IC_M [-2])= {\mathbb D}_{M_s} (IC_{M_s} [-1])=
({\mathbb D}_{M_s} (IC_{M_s})) [1];$
\item similar to (2), but for $\psi$
instead of $i^*.$
\end{enumerate}  
We thus get that $\iota_M$
induces an isomorphism of pairings
between $({\iota_M})|{M_s}$
and $(\iota_M)|{M_{\overline{\eta}}}$ via $i^* \to \psi.$
These pairing are the respective duality  pairings
$\iota_{M_s}$ and $\iota_{M_{\overline{\eta}}}$.
Indeed, the duality pairing isomorphisms at the level of intersection complexes are uniquely determined by their restriction to the smooth dense open subsets, where we see that they all preserve the unit section.

To conclude the statement for any two geometric fibers, we can find specializations to closed points, choose DVRs realizing the specializations, and repeat the argument above.
\end{proof}

\subsection{\'Etale cohomology of families of good moduli spaces}
For the following result, we work with $\overline{\mathbb{Q}}_{\ell}$-sheaves for the \'etale topology.
\begin{thm} \label{thm: coh of gms in families}
Let $f: \cX \to \cS$ be a smooth morphism between finite type algebraic stacks with affine diagonal over $k$. Suppose that $\cX$ admits a relative good moduli space $\overline{f}: \cM \to \cS$ such that $\overline{f}$ is proper and schematic. Suppose that $\cS$ is connected. Then the cohomology rings of any two geometric fibers  of $\overline{f}$
are non-canonically isomorphic. Their intersection
cohomology groups, endowed with their natural  module over the respective cohomology rings structures, and with their natural intersection pairings, are also non-canonically isomorphic.

If, furthermore, $\cS$ is smooth, then the following hold.

\begin{enumerate}
\item[(a)]
For every $i,$ the direct image sheaves
$R^i \overline{f}_* (\overline{\mathbb{Q}}_{\ell})_{\cM}$ and
$R^i \overline{f}_* \mathcal{IC}_{\cM}$
are lisse on $\cS$.

\item[(b)]
We have the identity of perverse sheaves on $\cS$
\[
(R^i\overline{f}_* (\mathcal{IC}_{\cM}))[dim(\cS)] =
{}^p\mathcal{H}^{dim(\cS)+i} (R\overline{f}_* (\mathcal{IC}_{\cM})).
\]
\item[(c)]
We have an isomorphism
\[
R\overline{f}_* \mathcal{IC}_{\cM} \simeq \bigoplus_{i\geq 0}
R^i\overline{f}_* \mathcal{IC}_{\cM}[-i].
\]
\end{enumerate}
\end{thm}
\begin{proof}
    In view of \Cref{cor: equisingularity}, the morphism $\overline{f}: \cM \to \cS$ is locally isotrivial as in \Cref{defn: locally isotrivial}. Therefore, the theorem follows directly from \Cref{supp coeff} and \Cref{iso coh fibers}.
\end{proof}

 \subsection{Singular cohomology in families over \texorpdfstring{$\mathbb C$}{\textbf{C}} and topological local triviality} \label{subsection: classical}

 In this section the ground field  is $k=\mathbb C.$
The relevant categories for a complex algebraic variety $X$ are: 
\begin{enumerate}
    \item the $\overline{\mathbb Q}_\ell$-constructible derived category
$D^b_c(X, \overline{\bQ}_\ell)$ (\'etale topology);
\item the algebraically constructible bounded derived  categories $\widetilde{D^b_c} (-, \bQ)$ with rational (or other) coefficients (classical topology); 
\item  M. Saito's bounded derived category $D^b(\MHM_{\mathrm{alg}} (-, \bQ))$ of algebraic mixed Hodge modules with rational coefficients (cf. \cite{saito-mixed-hodge, schnell-overview-saito}).
\end{enumerate}
They are related as follows
\begin{equation}\label{cats}
\xymatrix{
D^b(\MHM_{\mathrm{alg}} (X, \bQ))
\ar[r]^-{\operatorname{rat}} &
\widetilde{D^b_c}(X, \bQ) 
\ar[r]^-{
\otimes_{\bQ} \overline{\bQ}_\ell}
&
\widetilde{D^b_c}(X, \overline{\bQ}_\ell)
&
D^b_c (X, \overline{\bQ}_\ell) 
\ar[l]_-{\epsilon ^*}.
}
\end{equation}
Here, $\epsilon^*$ is the Artin comparison
fully faithful functor in \cite[\S 6.1.2]{bbdg}; the symbol $\otimes_{\bQ} \overline{\bQ}_\ell$ denotes the change-of-coefficients functor; and $\operatorname{rat}$ is the faithful functor
in the main result in \cite[Thm. 0.1]{saito-mixed-hodge}. 

These categories are equipped with their natural nine functors: $f^*$, $Rf_*$,
$Rf_!$, $f^!$, $R{\cH}om$, $\otimes^{\bL}$, the Verdier Duality functor, and the vanishing and nearby functors.  The functors in
\eqref{cats} are compatible (commute) with these  nine functors. Note that the following hold:
\begin{enumerate}[(1)]
    \item $\epsilon^*L$ is  locally constant sheaf if and only if  $L$ is lisse; 
    \item there is a direct sum decomposition $u: \epsilon^* A \simeq B\oplus C$ if and only if there is a
direct sum decomposition $u': A \simeq B' \oplus C'$ with $u=\epsilon^*(u');$
\item taking direct images commutes with changing coefficients; 
\item a sheaf of $\bQ$-vector spaces is locally constant  if and only if it is locally constant after changing coefficients
to $\overline{\bQ}_{\ell};$
\item changing coefficients and $\epsilon^*$ preserve the constant sheaf and the intersection complex.
\end{enumerate}
In particular, 
if a direct image sheaf of the constant sheaf or of the intersection complex  is lisse,
then the corresponding direct image is locally constant.

The aforementioned nine functors on the constructible derived categories $\widetilde{D}^b(-,  \bQ )$ admit canonical lifts along  the functor ${\rm rat}$
to M. Saito's categories $D^b ({\rm MHM}_{\mathrm{alg}}(-, \bQ)$. If we denote 
the lift of any of these functors $F$ by $\underline F$, then we have identities of the form $\operatorname{rat} \circ \underline{F} = F \circ \operatorname{rat}$. For example, for a morphism
$f:X \to Y$ of complex algebraic varieties, we denote by $\underline{R f_*}:
D^b (\MHM_{\mathrm{alg}} (X, \bQ )) \to D^b (\MHM_{\mathrm{alg}} (Y, \bQ ))  $  the canonical lift of the functor $Rf_*:
D^b_c (X, \bQ) \to D^b_c(Y, \bQ)$;
then we have that $Rf_* (\operatorname{rat}(A)) = \operatorname{rat}(\underline{Rf_*}(A))$ for every complex of mixed Hodge modules $A$ on $X$.

By construction, the standard $t$-structure on $D^b(\MHM_{\mathrm{alg}} (X, \bQ))$ corresponds to the (middle)-perversity $t$-structure on
$D^b_c(X, \bQ)$. 
In particular, given a complex $A$ of mixed modules on $X$, a morphism $f:X \to Y$ of complex algebraic varieties,  we have that
$^p\! \cH^i (Rf_* (\operatorname{rat} (A)) ) = 
\operatorname{rat} (\cH^{i} (\underline{Rf_*} A))$. 
The functors in \eqref{cats} are $t$-exact
for the corresponding $t$-structures
(standard for mixed Hodge modules, perverse for the others).

For an arbitrary complex algebraic variety $X$,
the sheaf $\bQ_X$ and the intersection complex $\cI \cC_X$
admit canonical lifts to $D^b(\MHM_{\mathrm{alg}} (X, \bQ))$
denoted by $\underline{\bQ_X}$ and 
$\underline{\cI \cC_X}$.

{\textbf{Smooth mixed Hodge modules.}} A mixed Hodge module $A$ on a complex variety $X$  is smooth  on $X$ if $L:=\mathrm{rat}(A)[-\dim (X)]$
is locally constant on $X$ (concentrated in cohomological degree zero)
and $X$ is smooth of pure dimension. 
In this case, $L$ is an admissible variation of mixed Hodge structures  in the sense of Zucker and of Kashiwara; the converse holds (cf. \cite[Thm. 0.2 and the paragraph preceding it]{saito-mixed-hodge}).
The constant sheaf-type object $\underline{\bQ_X}$ is a complex of mixed Hodge modules of weights at most $0$ (cf.
\cite[\S4.5]{saito-mixed-hodge}). Having weights at most $w$ is preserved by  pushforward via a proper morphism \cite[\S4.5.2]{saito-mixed-hodge}.

 \textbf{Smooth polarizable pure Hodge modules.}
A polarizable pure  Hodge module $A$ of weight $w$ on a complex variety $X$  is smooth  on $X$ if $L:=\mathrm{rat}(A)[-\dim (X)]$
is lisse on $X$ (concentrated in cohomological degree zero)
and $X$ is smooth of pure dimension.
In this case, $L$ is a polarizable variation of pure Hodge structures of weight $w-\dim (X)$; the converse holds  (cf. \cite[Lm. 3 and Thm. 2]{saito-mhp}). 
The intersection complex-type  object $\underline{{\cI  \cC}_X} [\dim (X)]$ 
is a polarizable pure Hodge module of weight $\dim (X)$. 

 \textbf{M. Saito's Decomposition Theorem.}
Bounded complexes of polarizable pure Hodge modules are preserved via pushforward by a projective  morphism
(cf. \cite[Thm. 1]{saito-mhp}) so that the
M. Saito's  Decomposition Theorem holds: if $A$ is a complex of 
polarizable pure Hodge modules on a complex variety $X$ and $f:X\to Y$
is a morphism of complex algebraic varieties, then $\underline{Rf_*} A
\simeq \oplus_i \cH^i (\underline{Rf_*} A) [-i]$ splits into the direct sum of shifted polarizable pure Hodge modules.

\textbf{Cohomology mixed Hodge modules/sheaves and smooth mixed Hodge modules.}  If
$X$ is smooth of pure dimension, $A$ is a smooth complex of mixed Hodge modules, i.e. $\mathrm{rat}(A)$
has lisse cohomology sheaves $\cH^i(\mathrm{rat}(A))$ on $X$, then we have the identity of perverse sheaves
$^p\!\cH^{\dim (X) +i} (\mathrm{rat}(A))= \cH^j(\mathrm{rat}(A)) [\dim (X)]$, so that the cohomology mixed Hodge modules $\cH^{\dim (X) +i}(A)$
is smooth on $X$ and the lisse sheaf
$\cH^i(\mathrm{rat}(A))$ underlies an admissible variation of mixed Hodge structures (pure,
if $A$ is a complex f polarizable pure Hodge modules).

\begin{prop} \label{prop: mixed hodge module prop}
   Let $f: M \to S$ be a locally isotrivial morphism (\Cref{defn: locally isotrivial}). Assume that $M$ and $S$ are complex algebraic varieties, that $f$ is projective, and that $S$ is smooth and connected.
The evident mixed Hodge modules analogue of \Cref{supp coeff} holds, namely we have the following.

\begin{enumerate}
\item[(a)]
For every $i \in \bZ^{\geq 0}$, the direct image sheaves
$R^i f_* (\overline{\mathbb{Q}}_{\ell})_{M}$ and
$R^i f_* \mathcal{IC}_{M}$
are locally constant on $S$. 
For each $i$, the sheaf $R^i f_* (\overline{\mathbb{Q}}_{\ell})_{M}$ underlies
a smooth mixed Hodge module $\underline{L[\dim (S)]}$ on $S$, so that the local system $L$ carries a  natural  admissible mixed Hodge structure 
of weights at most $i$; similarly, for the
sheaf $R^i f_* \mathcal{IC}_{M}$,
except that in this case it underlies a polarizable variation of pure Hodge structures of weight $i$.

\item[(b)]
We have the identity of perverse sheaves on $S$
\[
(R^if_* (\mathcal{IC}_{M}))[{\rm dim} (S)] =
{}^p\mathcal{H}^{\dim (S)+i} (Rf_* (\mathcal{IC}_{M})),
\]
as well as the identity
in the category 
$\MHM_{\mathrm{alg}}(S,\bQ)$
\[
\underline{(R^if_* (\mathcal{IC}_{M}))}[{\rm dim} (S)] =
\cH^{\dim (S)+i} (\underline{Rf_*} \, (\underline{\mathcal{IC}_{M}})).
\] 
\item[(c)]
We have an isomorphism
in $\widetilde{D^b_c}(S, \bQ)$
\[
Rf_* (\mathcal{IC}_{M}) \simeq \bigoplus_{i\geq 0}
R^if_* (\mathcal{IC}_{M})[-i],
\]
as well as the corresponding isomorphism in $D^b(\MHM_{\mathrm{alg}}(S, \bQ))$
\[
\underline{Rf_*} \, 
(\underline{{\cI \cC}_{M}}) \simeq \bigoplus_{i\geq 0}
\cH^{\dim (S) +i} (\underline{Rf_*} \, (\underline{\mathcal{IC}_{M}}))[-i].
\]
\end{enumerate}
\end{prop}

\begin{proof}
\Cref{supp coeff}, implies, via application of the comparison functor $\epsilon^*$,
 the following  statements
   at the level of
 $\widetilde{D^b_c}(S, \overline{\mathbb Q}_\ell)$, and thus at the level of 
 $\widetilde{D^b_c}(S,  \bQ)$:
 the local constancy assertions in (a); the identity of perverse sheaves in (b); the isomorphism  in $\widetilde{D^b_c}(S, \bQ)$
 appearing in (c).
 
Part (a) follows from the characterizations of
smooth mixed Hodge modules and of smooth polarizable pure Hodge modules.
Part (b) follows from the structure of the cohomology of a smooth  complex of mixed Hodge modules.
Part (c) follows from part (b) and  from  M. Saito's Decomposition Theorem.
\end{proof}

Next, we prove the following  topological local triviality result.

\begin{prop}[Local triviality of families]\label{loc triv}
Let $f: M \to S$ be a locally projective morphism of $\mathbb{C}$ varieties. Suppose that $f$ is locally isotrivial (\Cref{defn: locally isotrivial}) and that $S$ is smooth. The family $f: M \to S$
is topologically locally trivial over the base $S$.
\end{prop}
\begin{proof} We work with the classical topology.
Without loss of generality, we may assume that the morphism
$f$ factors as a closed embedding into $\mathbb{P}^n_S$,
followed by the restriction of the projection onto $S.$
The desired conclusion is a consequence of 
the First Thom Isotopy Lemma \cite[pg. 41, \S1.5, Thm.]{stratified-morse} if we can produce a Whitney stratification  $\mathfrak S$ of $M$ all of whose strata $T$ map submersively onto $S.$

It is a deep result of B. Teissier
\cite{teissier-polar} (see \cite[pg. 8]{decataldo-migliorini-projectors} for a short discussion)
that every complex  variety  admits a canonical algebraic Whitney stratification,
which we call the WT-stratification.
In order to prove the proposition, it is enough  that the strata of the WT-stratification of $M$ map submersively onto $S.$ We do so by using diagram \eqref{diagram: local isotriviality}, which we may arrange to be a diagram of varieties $T,F,U,V$ in this case. 

The WT-stratification is preserved via
local complex analytic isomorphisms.
It follows that: if $U\to M$ is \'etale, then the pre-image of 
the WT-stratification of $M$ is the
WT-stratification $U;$ if $U\to M$ is an \'etale cover, then a stratification of $M$ is the WT-stratification of $M$ if and only if its pre-image in $U$ is the WT-stratification of $U$;
the WT-stratification of a product $T\times F$, with $T$ smooth,   is the pre-image of the WT-stratification of $F$ via the projection. 

It is clear that the WT-stratification
of $T\times F$
has strata that map submersively onto $T.$ 
In view of \Cref{cor: equisingularity}, the discussion in the previous paragraph implies that the  WT-stratification of $M$ has strata that map submersively onto $S.$ 
\end{proof}

\begin{notn}
    Let $M$ be a finite type scheme over $\mathbb{C}$. We write $H^{\bullet}(M, \mathbb{Q})$ and $I\!H^{\bullet}(M, \mathbb{Q})$ for the corresponding  singular (intersection) cohomology groups with rational coefficients. We shall omit the coefficients $\mathbb{Q}$ from the notation whenever it is clear from context.
\end{notn} 

We end this section with the following result, where we denote by $h^{p,q}(\mathrm{Gr}^W_a H)$
the dimension of the $(p,q)$-part of
the weight $a$-graded piece of a rational
mixed Hodge structure $H$.
\begin{thm}\label{thm: two fibres of locally isotrivial family are cohomologically equivalent singular cohomology}
    Let $f: \cM \to \cS$ be a schematic locally projective morphism of quasi-separated stacks of finite type over $\mathbb{C}$. Suppose that $f$ is locally isotrivial and that $\cS$ is connected. Then, for any two $\mathbb{C}$-points $s,s' \in \cS(\mathbb{C})$, the following hold.
    \begin{enumerate}[(a)]
        \item The projective $\mathbb{C}$-varieties $\cM_s$ and $\cM_{s'}$ are homeomorphic.
        \item There is a noncanonical isomorphism $H^{\bullet}(\cM_s) \xrightarrow{\sim} H^{\bullet}(\cM_{s'})$ that preserves the weight filtrations. Such noncanonical isomorphism can be arranged to come from a homeomorphism as in (a).
        \item There is a noncanonical isomorphism $I
        \!H^{\bullet}(\cM_s) \xrightarrow{\sim} I\!H^{\bullet}(\cM_{s'})$ that is compatible with the
        respective intersection pairings and with the 
        structures of modules over the respective cohomology rings intertwined by an isomorphism as in (b). 
        \item For any indexes $p,q$ we have 
        \[h^{p,q}(Gr^W_a H^{\bullet}(\cM_s))= h^{p,q}(Gr^W_a H^{\bullet}(\cM_{s'})), \; \; \; \text{and} \; \; \; h^{p,q}(I\!H^{\bullet}(\cM_s))= h^{p,q}(I\!H^{\bullet}(\cM_{s'})).\]
        In particular, the corresponding $E$
polynomials of the mixed Hodge structures $H^{\bullet}(-)$ and of the pure Hodge structures  $I\!H^{\bullet}(-)$ of the two fibers agree.

    \end{enumerate}
\end{thm}
\begin{proof}
After pulling back to a resolution of singularities $S' \to S$ of a connected component of an atlas $S \to \cS$, we may assume that $f: M \to S$ is a locally projective morphism of schemes with $S$ smooth. In that case, part (a) follows from \Cref{loc triv}. Our proof of the local constancy of the direct image sheaves
$R^i f_* \overline{\mathbb{Q}}_M$,
conjoined with \cite[Prop. 6.2.3]{bbdg}, implies that
the weight filtrations on the stalks give rise to a locally constant filtration of the locally constant
$R^i f_* \overline{\mathbb{Q}}_M$. This implies part (b). Part (c) follows from part (a). Finally, (d) is a consequence of \Cref{prop: mixed hodge module prop}.
\end{proof}


\section{Applications}\label{sec: applications}

In this section, we assume that the algebraically closed field $k$ has characteristic $0$ and we fix a prime number $\ell$. We will exhibit several examples of pairs of schemes $M_1, M_2$ of finite type over $k$ having the same cohomology rings and intersection cohomology groups.

For the sake of precision, we introduce the following property.

\begin{property}[For $M_1$ and $M_2$]\label{property: cohomologically equivalent l-adic} There exist:
    \begin{enumerate}[(1)]
    \item\label{item: isomorphism of l-adic cohomology rings} a $\overline\bQ_\ell$-algebra isomorphism $H^{\bullet}\left(M_1, \overline{\mathbb{Q}}_{\ell}\right) \cong H^{\bullet}\left(M_2, \overline{\mathbb{Q}}_{\ell}\right)$ between $\ell$-adic cohomology rings;
        
    \item a $H^{\bullet}\left(M_1, \overline{\mathbb{Q}}_{\ell}\right)$-module isomorphism $I\!H^{\bullet}\left(M_1, \overline{\mathbb{Q}}_{\ell}\right) \cong I\!H^{\bullet}\left(M_2, \overline{\mathbb{Q}}_{\ell}\right)$, compatible with the intersection pairing in intersection cohomology.
    \end{enumerate}

    Here, the $H^{\bullet}\left(M_1, \overline{\mathbb{Q}}_{\ell}\right)$-module structure on $I\!H^{\bullet}\left(M_2, \overline{\mathbb{Q}}_{\ell}\right)$ is induced from its natural $H^{\bullet}\left(M_2, \overline{\mathbb{Q}}_{\ell}\right)$-module structure by transport along the isomorphism of cohomology rings in \eqref{item: isomorphism of l-adic cohomology rings}.
\end{property}

When working over the complex numbers, $k=\bC$, we consider the analogous property for singular cohomology, intersection cohomology and (mixed) Hodge structures on them.

\begin{property}[For $M_1$ and $M_2$, assuming $k=\bC$]\label{property: cohomologically equivalent over C}
There exist
  \begin{enumerate}[(1)]
        \item a homeomorphism $M_1^{\text{an}}\cong M_2^{\text{an}}$ between the underlying topological spaces of the analytifications of $M_1$ and $M_2$;
        
        \item \label{item: isomorphism singular cohomology rings} an isomorphism of rational cohomology rings \[H^{\bullet}\left(M_1, \mathbb{Q}\right) \cong H^{\bullet}\left(M_2, \mathbb{Q}\right)\] induced by the homeomorphism $M_1^{\text{an}}\cong M_2^{\text{an}}$ above, 
        that preserves the weight filtrations
        for the corresponding mixed Hodge structures, and

        \item an $H^{\bullet}\left(M_1, \mathbb{Q}\right)$-module isomorphism \[I\!H^{\bullet}\left(M_1, \mathbb{Q}\right) \cong I\!H^{\bullet}\left(M_2, \mathbb{Q}\right)\] between rational intersection cohomology groups, where the module structure on $I\!H^{\bullet}(M_2, \mathbb{Q})$ is induced by the isomorphism of cohomology rings in \eqref{item: isomorphism singular cohomology rings}.
    \end{enumerate}

    Moreover, the dimensions $h^{p,q} \left(Gr^W_a H^j\left(M_1,\bQ\right)\right)$ and $h^{p,q} \left( Gr^W_a H^j\left(M_2,\bQ\right)\right)$ coincide; in particular, the corresponding E-polynomial coincide. Likewise, the dimensions $h^{p,q}\left(I\!H^{\bullet}\left(M_1, \mathbb{Q}\right)\right)$ and $h^{p,q}\left(I\!H^{\bullet}(M_2, \mathbb{Q})\right)$ coincide.
\end{property}

\subsection{Moduli space of \texorpdfstring{$G$}{G}-bundles} \label{section: application G-bundles}

Fix a connected reductive group $G$ over $k$. Let $C$ be a smooth connected projective curve over $k$. The moduli stack $\Bun_G(C)$ of $G$-bundles on $C$ is a smooth algebraic stack with affine diagonal over $k$ (\cite[Prop. 1]{heinloth-uniformization}). The connected components of $\Bun_G(C)$ are in correspondence with elements of the algebraic fundamental group $\pi_1(G)$ \cite[Thm. 5.8]{hoffmann-connected-components}; for any $d \in \pi_1(G)$ we denote by $\Bun_G^d(C)$ the corresponding component of $G$-bundles of degree $d$. There is a notion of semistability for $G$-bundles defined by Ramanathan \cite{ramanathan-stable} which generalizes semistability of vector bundles on $C$. The locus of semistable $G$-bundles in $\Bun_G(C)$ is open, and corresponds to an open substack $\Bun_{G}^{d}(C)^{ss} \subset \Bun_G^d(C)$. The stack $\Bun_G^d(C)^{ss}$ admits a projective good moduli space over $k$ \cite{ramanathan-thesis-i, ramanathan-thesis-ii}, which we  denote by $M_G^d(C)$.

\begin{thm} \label{thm: moduli of g-bundles general k}
    Let $C_1, C_2$ be two smooth projective curves of the same genus over an algebraically closed field $k$ of characteristic $0$. Let $G$ be a connected reductive group over $k$, and fix $d \in \pi_1(G)$. Then, the moduli spaces $M^d_G(C_1)$ and $M^d_G(C_2)$ of semistable $G$-bundles on $C_1$ and $C_2$ satisfy \Cref{property: cohomologically equivalent l-adic}. If $k=\bC$, then $M^d_G(C_1)$ and $M^d_G(C_2)$ satisfy $\Cref{property: cohomologically equivalent over C}$.
\end{thm}
\begin{proof}
    Let $g$ denote the genus of the curves $C_1$ and $C_2$. 
    The stack $\cM_g$ of smooth curves of genus $g$ is a smooth, separated and irreducible Deligne--Mumford stack of finite type over $k$ \cite{deligne_mumford}. 
    Let $U\to \cM_g$ be an étale cover by a quasi-compact scheme $U$. 
    Since $\cM_g$ is connected, we may assume that the points of $\cM_g$ corresponding to $C_1$ and $C_2$ lie in the image of a connected component $S$ of $U$.
    
    Let $C\to S$ be the family of curves given by the map $S\to \cM_g$. 
    There is an algebraic stack $\Bun_G(C/S) \to S$ parametrizing $G$-bundles on the fibers of $C \to S$; this stack is locally of finite type and has affine relative diagonal over $S$ (by \cite[Theorem 1.3]{hall-rydh-tannakahom}, with $S=A$, $Z=C \times A$ and $X=BG$). This stack is known to be smooth (see the proof of \cite[Prop. 1]{heinloth-uniformization}, which applies over a general base scheme $S$). There is an open and closed substack $\Bun_G^d(C/S) \subset \Bun_G(C/S)$ parametrizing $G$-bundles that have degree $d$ in the sense of \cite{hoffmann-connected-components}. It is known that the subfunctor $\Bun_G^d(C/S)^{ss}\subset \Bun_G^d(C/S)$ parametrizing semistable $G$-bundles of degree $d$ is represented by an open substack (see for example \cite{gurjar-nitsure-higher} for opennes of semistability in the setting of families). It follows as a special case of \cite[Thm. 1.1]{langer-moduli-bundles-singular} that $\Bun_G^d(C/S)^{ss}$ admits a locally projective good moduli space $ f: M_G^d(C/S) \to S$ (alternatively, the GIT constructions of Schmitt \cite{schmitt-decorated-book} can be modified to apply to over such base $S$ using Seshadri's relative GIT \cite{seshadri-relative-git}).
    
    Let $s_1, s_2 \in S(k)$ be two points such that the fibers of $C\to S$ over $s_1$ and $s_2$ are isomorphic to $C_1$ and $C_2$, respectively.
    Since the formation of the good moduli space of $\Bun_G(C/S)_d^{ss}$ commutes with base-change on the base $S$ \cite[Prop. 4.7(i)]{alper-good-moduli}, the fibers of $f$ over $s_1$ and $s_2$ are isomorphic to $M_{G}^d(C_1)$ and $M_{G}^d(C_2)$, respectively. By applying \Cref{iso coh fibers} to $f$, we get that $M_{G}^d(C_1)$ and $M_{G}^d(C_2)$ satisfy \Cref{property: cohomologically equivalent l-adic}.

    If $k=\bC$ then, by \Cref{thm: two fibres of locally isotrivial family are cohomologically equivalent singular cohomology}, $M_{G}^d(C_1)$ and $M_{G}^d(C_2)$ satisfy \Cref{property: cohomologically equivalent over C}
\end{proof}

\begin{remark}\label{okn}
It follows as in \Cref{loc triv} that the moduli spaces of semistable $G$-bundles on any two curves of the same genus are homeomorphic. This can also be seen directly as a consequence of the Narasimhan--Seshadri correspondence for $G$-bundles (see e.g. \cite[Thm. 1.0.3(3)]{balaji_seshadri}).
\end{remark}

\subsection{Moduli of sheaves on surfaces with \texorpdfstring{$\omega$}{ω}-negative polarization} \label{application: negatively polarized surfaces}

In this subsection, we show that the (intersection) cohomology and the homeomorphism type of certain moduli spaces of sheaves on smooth surfaces vary nicely in families. Our key assumption will be that the polarization used to define semistability of sheaves satisfies the following property.
\begin{defn}
Let $X$ be a smooth projective surface over a field $K$ equipped with an ample line bundle $\cL$. We say that $\cL$ is an $\omega$-negative polarization if $\cL \cdot \omega_X <0$, where $\omega_X = \text{det}(\Omega^1_{X/K})$ denotes the canonical bundle of $X$. In this case, we say that the pair $(X, \cL)$ is an $\omega$-negatively polarized surface.
\end{defn}

The main motivation for this definition is that it ensures that $\Ext^2_{X}(\mathcal{F}, \mathcal{F})=0$ for any Gieseker $\mathcal{L}$-semistable sheaf of dimension 2 or 1. See the proof of \Cref{prop: smoothness of moduli of torsion-free sheaves for negatively polarized surfaces}.

\begin{example}
The following are examples of $\omega$-negatively polarized surfaces.
\begin{enumerate}[(a)]
    \item Any smooth del Pezzo surface $X$ admits the $\omega$-negative polarization $\omega_{X}^{-1}$.
    \item All iterated blowups of $\mathbb{P}^2$ at finitely many points admit $\omega$-negative polarizations.
    \item Any ruled smooth surface admits an $\omega$-negative polarization.
\end{enumerate}
\end{example}

\begin{remark}
    The Enriques--Kodaira classification of surfaces implies that a smooth surface over $k$ admits a negative polarization if and only if its Kodaira dimension is $-\infty$. Therefore the examples (b) and (c) above are exhaustive.
\end{remark}

\begin{defn}
    Let $S$ be a Noetherian scheme over $k$. A family of $\omega$-negatively polarized surfaces consists of a pair $(X \to S, \cL)$ where $X \to S$ is a smooth proper family equipped with a relatively ample line bundle $\cL \in \Pic(X)$ such that for all $s \in S$ the fiber $(X_s, \cL_s)$ is an $\omega$-negatively polarized surface.
\end{defn}

Let $(X \to S, \cL)$ be a family of $\omega$-negatively polarized surfaces over a Noetherian $k$-scheme $S$, and fix a polynomial $P \in \mathbb{Q}[n]$. We denote by $\Coh^{P, ss}_{X/S}$ the stack of Gieseker $\cL$-semistable pure sheaves on $X$ with Hilbert polynomial $P$. Recall that a coherent sheaf $\cF$ on a variety over a field is called pure if for all coherent subsheaves $0 \neq \cE \subset \cF$ we have $\text{dim}(\text{Supp}(\cE)) = \text{dim}(\text{Supp}(\cF))$.

To be more precise, for any $S$-scheme $T \to S$ the groupoid $\Coh^{P, ss}_{X/S}$ consists of $T$-flat families of finitely presented quasicoherent $\cO_{X_T}$-modules $\cF$ on $X_T$ such that for all geometric points $\overline{t} \to T$ the base-change $\cF|_{X_{\overline{t}}}$ is a Gieseker $\cL_{X_{\overline{t}}}$-semistable pure sheaves with Hilbert polynomial $P$. The stack $\Coh^{P, ss}_{X/S} \to S$ has affine relative diagonal and is of finite type over $S$, as it is an open substack of the stack of coherent sheaves which satisfies both properties \cite[\href{https://stacks.math.columbia.edu/tag/08KA}{Tag 08KA}]{stacks-project}. If the field $k$ is of characteristic $0$, then $\Coh^{P, ss}_{X/S}$ admits a good moduli space $M^P_{X/S}$ which is locally projective over $S$ (\cite[\S1]{Simpson-repnI}, \cite[Thm. 1.1]{langer-lie-algebroids}, \cite[Thm. 5.2]{torsion-freepaper}).

Note that, if the Hilbert polynomial $P$ has degree 2, then $\Coh^{P, ss}_{X/S}$ parametrizes torsion-free sheaves.

\begin{prop} \label{prop: smoothness of moduli of torsion-free sheaves for negatively polarized surfaces}
    Let $(X \to S, \cL)$ be a family of $\omega$-negatively polarized surfaces over a Noetherian $k$-scheme $S$. Suppose that the Hilbert polynomial $P$ has degree 2. Then the stack $\Coh^{P, ss}_{X/S}$ is smooth over $S$.
\end{prop}
\begin{proof}
    By \cite[\href{https://stacks.math.columbia.edu/tag/0DP0}{Tag 0DP0} and \href{https://stacks.math.columbia.edu/tag/02HT}{Tag 02HT}]{stacks-project} applied to the finite type morphism $\Coh^{P, ss}_{X/S} \to S$, it suffices to show the infinitesimal lifting criterion for square-zero thickenings for small extensions  $A \twoheadrightarrow B$ of local Artin $k$-algebras over $S$. We denote by $\kappa$ the common residue field of $A$ and $B$, and we set $I$ to be the ideal of $B$ in $A$. Choose a morphism $\Spec(B) \to \Coh^{P, ss}_{X/S}$ corresponding to a $B$-flat family of sheaves $\cF$ on the base-change $X_B$. We denote by $\cF_{\kappa}$ the restriction to $X_{\kappa}$. By \cite[\href{https://stacks.math.columbia.edu/tag/08VW}{Tag 08VW}]{stacks-project}, there is an obstruction $\text{ob} \in \Ext^2_{X_{\kappa}}(\cF_{\kappa}, \cF_{\kappa}) \otimes_{\kappa} I$ such that $\text{ob} =0$ if and only if a lift $\Spec(A) \to \Coh^{P, ss}_{X/S}$ exists. We conclude the proof by showing that $\Ext^2_{X_{\kappa}}(\cF_{\kappa}, \cF_{\kappa}) =0$. By Serre duality for the smooth surface $X_{\kappa}$, this is isomorphic to $\Hom_{X_{\kappa}}(\cF_{\kappa}, \cF_{\kappa} \otimes \omega_{X_{\kappa}})$. By assumption $\cF_{\kappa}$ is a $\cL_{\kappa}$-Gieseker semistable sheaf on $X_{\kappa}$, and therefore the same holds for the twist $\cF_{\kappa} \otimes \omega_{X_{\kappa}}$. Since $\cL_{\kappa} \cdot \omega_{X_{\kappa}} <0$, it follows that the slope of $\cF_{\kappa} \otimes \omega_{X_{\kappa}}$ is strictly smaller than the slope of $\cF_{\kappa}$. Since both are Gieseker semistable, it follows that all homomorphisms from $\cF_{\kappa}$ to $\cF_{\kappa} \otimes \omega_{\kappa}$ are zero by \cite[Lemma 1.3.3]{huybrecths-lehn}.
\end{proof}

\begin{thm} \label{thm: moduli of negatively polarized surfaces}
    Let $(X \to S, \cL)$ be a family of $\omega$-negatively polarized surfaces over a connected finite type $k$-scheme $S$, where $k$ is an algebraically closed field of characteristic $0$. Fix a Hilbert polynomial $P$ of degree 2. For any two $k$-points $s, t \in S(k)$, let $M_{X_s}^{P}$ (resp. $M_{X_t}^{P}$) denote the corresponding moduli space of Gieseker $\cL_s$-semistable torsion-free sheaves on $X_s$ (resp. Gieseker $\cL_t$-semistable sheaves on $X_t$) \cite[\S1]{Simpson-repnI}. Then $M_{X_s}^{P}$ and $M_{X_t}^{P}$ satisfy \Cref{property: cohomologically equivalent l-adic}. If $k=\bC$, then $M_{X_s}^{P}$ and $M_{X_t}^{P}$ satisfy \Cref{property: cohomologically equivalent over C}.
\end{thm}
\begin{proof}
 By smoothness of the stack $\Coh^{P, ss}_{X/S}$, proven in \Cref{prop: smoothness of moduli of torsion-free sheaves for negatively polarized surfaces}, its good moduli space $M^P_{X/S}\to S$ is locally isotrivial in the sense of \Cref{defn: locally isotrivial} (\Cref{cor: equisingularity}). The fibers of $M^P_{X/S}\to S$ over $s$ and $t$ are $M_{X_s}^{P}$ and $M_{X_t}^{P}$ respectively, since the formation of good moduli spaces commutes with base change from $S$ \cite[Prop. 4.7(i)]{alper-good-moduli}. The result follows from \Cref{iso coh fibers} and, in the case $k=\bC$, from \Cref{thm: two fibres of locally isotrivial family are cohomologically equivalent singular cohomology}.
\end{proof}

By the smoothness of $X \to S$, there is a determinant morphism $\text{det}: \Coh_{X/S}^{P,ss} \to \Pic_{X/S}$ to the relative Picard scheme $\Pic_{X/S}$ \cite[\S2.1]{huybrecths-lehn}. The obstructions to lift sheaves $\cF_{\kappa}$ under this morphism belong to the trace-zero Ext group $\Ext_0^2(\cF_{\kappa}, \cF_{\kappa})$, and so in particular they are $0$ by the proof of \Cref{prop: smoothness of moduli of torsion-free sheaves for negatively polarized surfaces}. It follows that $\text{det}: \Coh_{X/S}^{P,ss} \to \Pic_{X/S}$ is smooth. By the universal property of good moduli spaces, there is an induced proper morphism $\overline{det}: M_{X/S}^P \to \Pic_{X/S}$. For a point $p\in \Pic_{X/S}(k)$ corresponding to a point $s\in S(k)$ and a line bundle $L$ on $X_s$, the fiber $\overline{det}^{-1}(p)$ is the moduli space $M_{X_s,L}^P$ of sheaves on $X_s$ with Hilbert polynomial $P$ and fixed determinant $L$.

\begin{prop}
In the situation of \Cref{thm: moduli of negatively polarized surfaces}, fix line bundles $L_1$ and $L_2$ on $X_s$ and $X_t$ respectively that are algebraically equivalent in the sense that the points of $\Pic_{X/S}$ that they lie in the same connected component. Then the moduli spaces $M_{X_s,L_1}^P$ and $M_{X_t,L_2}^P$ of sheaves with fixed determinant satisfy \Cref{property: cohomologically equivalent l-adic}. If $k=\bC$, then $M_{X_s,L_1}^P$ and $M_{X_t,L_2}^P$ satisfy \Cref{property: cohomologically equivalent over C}. 
\end{prop}
\begin{proof}
The map $\overline{det}: M_{X/S}^P \to \Pic_{X/S}$ is locally isotrivial by \Cref{cor: equisingularity}. Therefore we may apply \Cref{iso coh fibers} and, in the case $k=\bC$, \Cref{thm: two fibres of locally isotrivial family are cohomologically equivalent singular cohomology}, to conclude the result.
\end{proof}

 \textbf{Le Potier Morphism.}
Let $(X \to S, \cL)$ be a smooth family of $\omega$-negatively polarized surfaces, and fix a Hilbert polynomial $P$ of degree $1$. There is the Le Potier morphism $h: \Coh^{P, ss}_{X/S} \to \text{Hilb}^1_{X/S}$ originally defined by Le Potier \cite{le-potier-morphism}, where $\text{Hilb}^1(X/S)$ denotes the (relative) Hilbert scheme of subschemes of dimension 1 on the fibers of $X \to S$. For any $S$-scheme $T$, the morphism $h$ sends a $T$-flat family of pure sheaves of dimension 1 in $\Coh_{X/S}^{P, ss}(T)$ to its zero Fitting support, which in this case can be checked to be a relative Cartier divisor on $X_T \to T$ in $\text{Hilb}^1_{X/S}(T)$. By the universal property of good moduli spaces, this induces the proper morphism Le Potier
morphism $\overline{h}: M_{X/S}^P \to \text{Hilb}^1_{X/S}$. For every geometric point $s \to S$, the base change yields the  proper Le Potier morphism $\overline{h}_s: M_{X_s}^P \to \text{Hilb}^1_{X_s}$, which induces the perverse filtration on the intersection cohomology group (resp. cohomology group) $I\!H^{\bullet}\left(M_{X_s}^{P}, \overline{\mathbb{Q}}_{\ell}\right)$ (resp. $H^{\bullet}\left(M_{X_s}^{P}, \overline{\mathbb{Q}}_{\ell}\right)$). 

\begin{prop} \label{prop: le potier morphism filtration}
    Let $(X \to S, \cL)$ be a smooth family of $\omega$-negatively polarized over a connected finite type $k$-scheme $S$, where $k$ is an algebraically closed field of characteristic $0$. Fix a Hilbert polynomial $P$ of degree 1. Then, for any two points $s,t \in S(k)$, the corresponding moduli spaces $M_{X_s}^{P}$ and $M_{X_t}^{P}$ satisfy \Cref{property: cohomologically equivalent l-adic}, where in addition it can be ensured that the isomorphisms $IC^{\bullet}\left(M_{X_s}^{P},\overline{\bQ}_\ell\right)\cong IC^{\bullet}\left(M_{X_t}^{P}, \overline{\bQ}_\ell\right)$ and $H^{\bullet}\left(M_{X_s}^{P},\overline{\bQ}_\ell\right)\cong H^{\bullet}\left(M_{X_t}^{P}, \overline{\bQ}_\ell\right)$ are filtered isomorphisms  with respect to  the perverse filtrations induced by the corresponding Le Potier morphisms $\overline{h}_s$ and $\overline{h}_t$. 
    
    If $k=\bC$, then $M_{X_s}^{P}$ and $M_{X_t}^{P}$ satisfy \Cref{property: cohomologically equivalent over C} and the corresponding isomorphisms $IC^{\bullet}\left(M_{X_s}^{P}\right)\cong IC^{\bullet}\left(M_{X_t}^{P}\right)$ and $H^{\bullet}\left(M_{X_s}^{P}\right)\cong H^{\bullet}\left(M_{X_t}^{P}\right)$ can be arranged to be filtered isomorphisms with respect to  the perverse filtrations induced by the corresponding Le Potier morphisms.
\end{prop}

\begin{proof}
  By a similar reduction as in \Cref{iso coh fibers}, we can assume without loss of generality that $S$ is a smooth curve. Then, using the triviality of the vanishing cycle functor proven in \Cref{section: cohomology in families}, the proof of \cite[Thm. 3.2.1(i)]{de-cataldo-maulik-perverse} applies verbatim to show that the corresponding specialization isomorphisms from the (intersection) cohomology of any closed $S$-fiber to the (intersection) cohomology of the geometric generic fiber are compatible with the perverse filtrations: while the proof in \cite{de-cataldo-maulik-perverse} is in the context of the complex topology,  the argument found there is formal and also applies in the \'etale setting.
\end{proof}

\bibliographystyle{alpha}
\footnotesize{\bibliography{coh_G_bundles.bib}}

 \textsc{Department of Mathematics, Stony Brook University,
    Stony Brook, NY 11794-3651,
USA}\par\nopagebreak
  \textit{E-mail address}, \texttt{mark.decataldo@stonybrook.edu}
  
  \textsc{Department of Mathematics, University of Pennsylvania,
209 South 33rd Street,
Philadelphia, PA 19104, USA}\par\nopagebreak
  \textit{E-mail address}, \texttt{andresfh@sas.upenn.edu}

 \textsc{Mathematics Department, Columbia University, New York, NY 10027, USA}\par\nopagebreak
  \textit{E-mail address}, \texttt{andres.ibaneznunez@columbia.edu}
\end{document}